\documentclass[11pt]{article}%
\usepackage{amsmath}
\usepackage{amsfonts}
\usepackage{amssymb}
\usepackage{graphicx}
\usepackage{a4wide}
\usepackage{hyperref}%
\providecommand{\U}[1]{\protect\rule{.1in}{.1in}}
\newtheorem{theorem}{Theorem}

\newtheorem{corollary}[theorem]{Corollary}

\newtheorem{definition}[theorem]{Definition}
\newtheorem{example}[theorem]{Example}

\newtheorem{lemma}[theorem]{Lemma}

\newtheorem{proposition}[theorem]{Proposition}
\newtheorem{remark}[theorem]{Remark}

\newtheorem{question}{Question}

\begin{document}

\title{Series of convex functions: subdifferential, conjugate and applications to
entropy minimization}
\author{C. Vall\'{e}e\thanks{Institut PPRIME, UPR 3346, SP2MI, Futuroscope Chasseneuil
Cedex, France.}
\and C. Z\u{a}linescu\thanks{Faculty of Mathematics, University Al.\ I. Cuza,
700506 Ia\c{s}i, Romania, and Institute of Mathematics Octav Mayer (Romanian
Academy), Ia\c{s}i, Romania, e-mail: \texttt{zalinesc@uaic.ro}. The
contribution of C.Z. was done, mainly, during his 2014 visit at Institut
PPRIME, Universit\'{e} de Poitiers. }}
\date{}
\maketitle

\begin{abstract}
A formula for the sub\-differential of the sum of a series of convex
functions defined on a Banach space was provided by X. Y. Zheng in
1998. In this paper, besides a slight extension to locally convex
spaces of Zheng's results, we provide a formula for the conjugate of
a countable sum of convex functions. Then we use these results for
calculating the sub\-differentials and the conjugates in two
situations related to entropy minimization, and we study a concrete
example met in Statistical Physics.

\end{abstract}

\textbf{Key words:} Series of convex functions, sub\-differential,
conjugate, entropy minimization, statistical physics.

\section{Introduction}

The starting point of this study is the method used for deriving maximum
entropy of ideal gases in several books dedicated to statistical physics
(statistical mechanics); see \cite[pp.\ 119, 120]{LanLif:80}, \cite[pp.\ 15,
16]{Gue:07}, \cite[p.\ 43]{PatBea:11}, \cite[p.\ 39]{Bas:14}. The problem is
reduced to maximize $-\sum_{i\in I}n_{i}(\ln n_{i}-1)$ [equivalently to
minimize $\sum_{i\in I}n_{i}(\ln n_{i}-1)$] with the constraints $\sum_{i\in
I}n_{i}=N$ and $\sum_{i\in I}n_{i}\varepsilon_{i}=\varepsilon$ with $n_{i}$
nonnegative integers, or, by normalization (taking $p_{i}:=n_{i}/N$), to
maximize $-\sum_{i\in I}p_{i}(\ln p_{i}-1)$ [equivalently to minimize
$\sum_{i\in I}p_{i}(\ln p_{i}-1)$] with the constraints $\sum_{i\in I}p_{i}=1$
and $\sum_{i\in I}p_{i}e_{i}=e$ with $p_{i}\in\mathbb{R}_{+}$. For these one
uses the Lagrange multipliers method in a formal way. Even if nothing is said
about the set $I$, from examples (see \cite[(47.1)]{LanLif:80}, \cite[(3.11)]%
{Gue:07}, \cite[(1.4.5)]{PatBea:11}, etc) one guesses that $I$ is a countable
set. Our aim is to treat rigorously such problems. Note that the problem of
minimum entropy in the case in which the infinite sum is replaced by an
integral on a finite measure space and the constraints are defined by a finite
number of (continuous) linear equations is treated rigorously by J. M. Borwein
and his collaborators in several papers (in the last 25 years); see
\cite{Bor:12} for a recent survey.

The plan of the paper is the following. In Section 2 we present a slight
extension (to locally convex spaces) of the results of X. Y. Zheng
\cite{Zhe:98} related to the subdifferential of the sum of a series of convex
functions; we provide the proofs for readers convenience. In Section 3 we
apply the results in Section 2 for deriving a formula for the conjugate of the
sum of a series of convex functions, extending so Moreau's Theorem on the
conjugate of the sum of a finite family of convex functions to countable sums
of such functions. In Section 4 we apply the results in the preceding sections
to find the minimum entropy for a concrete situation from Statistical Physics.

\section{Series of convex functions}

Throughout this paper, having a sequence $(A_{n})_{n\geq1}$ of nonempty sets
and a sequence $(x_{n})_{n\geq1}$, the notation $(x_{n})_{n\geq1}\subset
(A_{n})_{n\geq1}$ means that $x_{n}\in A_{n}$ for every $n\geq1.$

In the sequel $(E,\tau)$ is a real separated locally convex space (lcs for
short) and $E^{\ast}$ is its topological dual. Moreover, we shall use standard
notations and results from convex analysis (see e.g.\ \cite{Mor:66},
\cite{Zal:02}). Consider $f_{n}\in\Lambda(E)$ (that is $f_{n}$ is proper and
convex) for every $n\geq1.$ Assume that $f(x):=\sum_{n\geq1}f_{n}%
(x):=\lim_{n\rightarrow\infty}\sum_{k=1}^{n}f_{k}(x)$ exists in $\overline
{\mathbb{R}}:=\mathbb{R}\cup\{-\infty,\infty\}$ for every $x\in E$, where
$\infty:=+\infty$. Then, clearly, the corresponding function $f:E\rightarrow
\overline{\mathbb{R}}$ is convex and $\operatorname*{dom}f\subset\cap_{n\geq
1}\operatorname*{dom}f_{n}$.

Note that, if $(f_{n})_{n\geq1}\subset\Gamma(E)$ (that is $f_{n}\in\Lambda(X)$
is also lower semi\-continuous, lsc for short) and there exists $(x_{n}^{\ast
})_{n\geq1}\subset(\operatorname*{dom}f_{n}^{\ast})_{n\geq1}$ such that the
series $\sum_{n\geq1}f_{n}^{\ast}(x_{n}^{\ast})$ is convergent and $w^{\ast}%
$-$\lim_{n\rightarrow\infty}\sum_{k=1}^{n}x_{k}^{\ast}=x^{\ast}\in E^{\ast}$
(that is $x^{\ast}=w^{\ast}$-$\sum_{n\geq1}x_{n}^{\ast}$), then $\lim
_{n\rightarrow\infty}\sum_{k=1}^{n}f_{k}(x)$ exists and belongs to
$(-\infty,\infty]$ for every $x\in E$; moreover, $f$ is lsc.

Indeed, since $g_{n}(x):=f_{n}(x)+f_{n}^{\ast}(x_{n}^{\ast})-\left\langle
x,x_{n}^{\ast}\right\rangle \geq0$, $g(x):=\lim_{n\rightarrow\infty}\sum
_{k=1}^{n}g_{n}(x)=\sup_{n\geq1}g_{n}(x)$ exists and belongs to $[0,\infty].$
But $\sum_{k=1}^{n}g_{n}(x)=\sum_{k=1}^{n}f_{n}(x)+\sum_{k=1}^{n}f_{n}^{\ast
}(x_{n}^{\ast})-\left\langle x,\sum_{k=1}^{n}x_{n}^{\ast}\right\rangle $, and
$\gamma:=\lim_{n\rightarrow\infty}\sum_{k=1}^{n}f_{n}^{\ast}(x_{n}^{\ast}%
)\in\mathbb{R}$, $\lim_{n\rightarrow\infty}\left\langle x,\sum_{k=1}^{n}%
x_{n}^{\ast}\right\rangle =\left\langle x,x^{\ast}\right\rangle $. It follows
that $f(x)=\lim_{n\rightarrow\infty}\sum_{k=1}^{n}f_{n}(x)=g(x)-\gamma
+\left\langle x,x^{\ast}\right\rangle \in(-\infty,\infty]$. Since $g_{n}%
\in\Gamma(X)$ for every $n$ and $g=\sup_{n\geq1}g_{n}$ is lsc, it follows that
$f$ is lsc, too.

\begin{definition}
\label{d-z1}\emph{(Zheng \cite[p.\ 79]{Zhe:98})} Let $A,A_{n}\in
\mathcal{P}_{0}(E):=\{F\subset E\mid F\neq\emptyset\}$ $(n\geq1)$. One says
that $(A_{n})_{n\geq1}$ \emph{converges normally} to $A$ (with respect to
$\tau$), written $A=\tau$-$\sum_{n\geq1}A_{n}$, if:

\begin{itemize}
\item[$(I)$] for every sequence $(x_{n})_{n\geq1}\subset(A_{n})_{n\geq1}$, the
series $\sum_{n\geq1}x_{n}$ $\tau$-converges and its sum $x$ belongs to $A;$

\item[$(II)$] for each ($\tau$-)neighborhood $U$ of $0$ in $E$ (that is
$U\in\mathcal{N}_{E}^{\tau}$) there is $n_{0}\geq1$ such that $\sum_{k\geq
n}x_{k}\in U$ for all sequences $(x_{n})_{n\geq1}\subset(A_{n})_{n\geq1}$ and
all $n\geq n_{0}$ (observe that the series $\sum_{k\geq n}x_{k}$ is $\tau
$-convergent by (I));

\item[$(III)$] for each $x\in A$ there exists $(x_{n})_{n\geq1}\subset
(A_{n})_{n\geq1}$ such that $x=\tau$-$\sum_{n\geq1}x_{n}$.
\end{itemize}
\end{definition}

Observe that $A$ in the above definition is unique; moreover, $A$ is convex if
all $A_{n}$ are convex.

\begin{remark}
\label{rem0}1) Assume that $E$ is the topological dual $X^{\ast}$ of the lcs
$X$ endowed with the weak$^{\ast}$ topology $w^{\ast}$, and $(A_{n})_{n\geq
1}\subset\mathcal{P}_{0}(X^{\ast})$ is such that (I) holds. Then (II) in
Definition \ref{d-z1} holds if and only if for every $\varepsilon>0$ and every
$x\in X$ there exists $n_{0}=n_{\varepsilon,x}\geq1$ such that $\big\vert\sum
_{k\geq n}\left\langle x,x_{k}^{\ast}\right\rangle \big\vert\leq\varepsilon$
for all sequences $(x_{n}^{\ast})_{n\geq1}\subset(A_{n})_{n\geq1}$ and all
$n\geq n_{0}.$

2) Assume that $E$ is a normed vector space (nvs for short) endowed with the
strong (norm) topology $s$, and $(A_{n})_{n\geq1}\subset\mathcal{P}_{0}(E)$ is
such that (I) holds. Then (II) in Definition \ref{d-z1} holds if and only if
for every $\varepsilon>0$ there exists $n_{\varepsilon}\geq1$ such that
$\big\Vert\sum_{k\geq n}x_{k}\big\Vert\leq\varepsilon$ for all sequences
$(x_{n})_{n\geq1}\subset(A_{n})_{n\geq1}$ and all $n\geq n_{\varepsilon}$. It
follows that $A=s$-$\sum_{n\geq1}A_{n}$ implies that $A$ is a
Hausdorff-Pompeiu limit of $(\sum_{k=1}^{n}A_{k})_{n\geq1}.$
\end{remark}

In the rest of this section we mainly reformulate the results of Zheng
\cite{Zhe:98} in the context of locally convex spaces without asking the
functions be lower semi\-continuous. We give the proofs for reader's convenience.

\begin{theorem}
\label{t-zt31}Let $f,f_{n}\in\Lambda(E)$ be such that $f(x)=\sum_{n\geq1}%
f_{n}(x)$ for every $x\in E$. Assume that $\overline{x}\in\operatorname*{core}%
(\operatorname*{dom}f)$. Then%
\[
f_{+}^{\prime}(\overline{x},u)=\sum_{n\geq1}f_{n+}^{\prime}(\overline
{x},u)\quad\forall u\in E.
\]

\end{theorem}

Proof. Consider $u\in E$. Because $\overline{x}\in\operatorname*{core}%
(\operatorname*{dom}f)$, there exists $\delta>0$ such that $\overline{x}%
+tu\in\operatorname*{dom}f\subset\operatorname*{dom}f_{n}$ for every $t\in
I:=[-\delta,\delta]$. Consider $\varphi,\varphi_{n}:I\rightarrow\mathbb{R}$
defined by $\varphi(t):=f(\overline{x}+tu)$, $\varphi_{n}(t):=f_{n}%
(\overline{x}+tu)$; $\varphi$, $\varphi_{n}$ are convex and $\varphi
(t)=\sum_{n\geq1}\varphi_{n}(t)$ for every $t\in I$. Of course, $f_{+}%
^{\prime}(\overline{x},u)=\lim_{t\rightarrow0+}\frac{\varphi(t)-\varphi(0)}%
{t}$, and similarly for $f_{n+}^{\prime}(\overline{x},u)$. Since the mappings
$I\setminus\{0\}\ni t\mapsto t^{-1}\varphi_{n}(t)\in\mathbb{R}$ are
nondecreasing we get
\[
\frac{\varphi_{n}(-\delta)-\varphi_{n}(0)}{-\delta}\leq\frac{\varphi
_{n}(t)-\varphi_{n}(0)}{t}\leq\frac{\varphi_{n}(\delta)-\varphi_{n}(0)}%
{\delta},
\]
whence%
\[
0\leq\psi_{n}(t):=\frac{\varphi_{n}(t)-\varphi_{n}(0)}{t}-\frac{\varphi
_{n}(-\delta)-\varphi_{n}(0)}{-\delta}\leq\frac{\varphi_{n}(\delta
)-\varphi_{n}(0)}{\delta}-\frac{\varphi_{n}(-\delta)-\varphi_{n}(0)}{-\delta
}=:\gamma_{n}%
\]
for all $n\geq1$ and $t\in(0,\delta]$. Since the series $\sum_{n\geq1}%
\gamma_{n}$ is convergent, the series $\sum_{n\geq1}\psi_{n}$ is uniformly
convergent on $(0,\delta]$. It follows that
\[
\lim_{t\rightarrow0+}\sum_{n\geq1}\psi_{n}(t)=\sum_{n\geq1}\lim_{t\rightarrow
0+}\psi_{n}(t).
\]
Since $\sum_{n\geq1}\frac{\varphi_{n}(-\delta)-\varphi_{n}(0)}{-\delta}%
=\frac{\varphi(-\delta)-\varphi(0)}{-\delta}$, we obtain that
\begin{align*}
f_{+}^{\prime}(\overline{x},u)  &  =\lim_{t\rightarrow0+}\frac{\varphi
(t)-\varphi(0)}{t}=\lim_{t\rightarrow0+}\sum_{n\geq1}\frac{\varphi
_{n}(t)-\varphi_{n}(0)}{t}=\sum_{n\geq1}\lim_{t\rightarrow0+}\frac{\varphi
_{n}(t)-\varphi_{n}(0)}{t}\\
&  =\sum_{n\geq1}f_{n+}^{\prime}(\overline{x},u).
\end{align*}
The proof is complete. \hfill$\square$

\begin{proposition}
\label{p-zl32}Let $f,f_{n}\in\Lambda(E)$ be such that $f(x)=\sum_{n\geq1}%
f_{n}(x)$ for every $x\in E$. Assume that the series $\sum_{n\geq1}f_{n}$
converges uniformly on a neighborhood of $x_{0}\in\operatorname*{int}%
(\operatorname*{dom}f)$. Then for every $x\in\operatorname*{int}%
(\operatorname*{dom}f)$ there exists a neighborhood of $x$ on which the series
$\sum_{n\geq1}f_{n}$ converges uniformly.
\end{proposition}

Proof. Replacing (if necessary) $f_{n}$ by $g_{n}$ defined by $g_{n}%
(x):=f_{n}(x_{0}+x)-f_{n}(x_{0})$ and $f$ by $g$ defined by $g(x):=f(x_{0}%
+x)-f(x_{0})$, we may (and do) assume that $x_{0}=0$ and $f_{n}(0)=f(0)=0.$
There exists a closed, convex and symmetric neighborhood $V$ of $x_{0}=0$ such
that $2V\subset\operatorname*{dom}f$ and the series $\sum_{n\geq1}f_{n}$
converges uniformly on $2V$. Set $p:=p_{V}$, the Minkowski functional
associated to $V$. Then $p$ is a continuous seminorm such that
$\operatorname*{int}V=\{x\in E\mid p(x)<1\}$ and $\operatorname*{cl}V=V=\{x\in
E\mid p(x)\leq1\}$. Consider $x\in\operatorname*{int}(\operatorname*{dom}f)$.
If $p(x)<2$ then $x\in\operatorname*{int}(2V)$; take $U:=2V$ in this case.

Let $p(x)\geq2$. Since $x\in\operatorname*{int}(\operatorname*{dom}f)$, there
exists $\mu>0$ such that $x^{\prime}:=(1+\mu)x\in\operatorname*{dom}f$, and so
$x=(1-\lambda)x^{\prime}+\lambda0$, where $\lambda:=\mu/(1+\mu)\in(0,1)$. Fix
$u\in V$ $(\Leftrightarrow p(u)\leq1)$; we have that $x+\lambda u=(1-\lambda
)x^{\prime}+\lambda u$, and so
\begin{equation}
f_{n}(x+\lambda u)\leq(1-\lambda)f_{n}(x^{\prime})+\lambda f_{n}%
(u)\quad\forall n\geq1. \label{r-z1}%
\end{equation}
On the other hand $1<2-\lambda\leq p(x)-p(\lambda u)\leq p(x+\lambda u)$, and
so $\frac{x+\lambda u}{p(x+\lambda u)}\in V$ and
\[
f_{n}\left(  \frac{x+\lambda u}{p(x+\lambda u)}\right)  \leq\frac
{1}{p(x+\lambda u)}f_{n}(x+\lambda u),
\]
whence
\begin{equation}
f_{n}(x+\lambda u)\geq p(x+\lambda u)f_{n}\left(  \frac{x+\lambda
u}{p(x+\lambda u)}\right)  \quad\forall n\geq1. \label{r-z2}%
\end{equation}
From (\ref{r-z1}) and (\ref{r-z2}) we get
\[
p(x+\lambda u)\sum_{k=l}^{m}f_{n}\left(  \frac{x+\lambda u}{p(x+\lambda
u)}\right)  \leq\sum_{k=l}^{m}f_{k}(x+\lambda u)\leq(1-\lambda)\sum_{k=l}%
^{m}f_{k}(x^{\prime})+\lambda\sum_{k=l}^{m}f_{k}(u),
\]
whence, because $p(x+\lambda u)\leq p(x)+1,$
\[
\bigg\vert\sum_{k=l}^{m}f_{k}(x+\lambda u)\bigg\vert\leq\bigg\vert\sum
_{k=l}^{m}f_{k}(x^{\prime})\bigg\vert+\bigg\vert\sum_{k=l}^{m}f_{k}%
(u)\bigg\vert+(p(x)+1)\bigg\vert\sum_{k=l}^{m}f_{n}\left(  \frac{x+\lambda
u}{p(x+\lambda u)}\right)  \bigg\vert
\]
for all $l,m\geq1$ with $l\leq m$. Since $u,\frac{x+\lambda u}{p(x+\lambda
u)}\in V\subset2V$, from the uniform convergence of $\sum_{n\geq1}f_{n}$ on
$2V$ and the convergence of $\sum_{n\geq1}f_{n}(x^{\prime})$, the previous
inequality shows that $\sum_{n\geq1}f_{n}$ is uniformly convergent on
$U:=x+\lambda V$ $(\subset\operatorname*{dom}f)$. The proof is complete.
\hfill$\square$

\begin{theorem}
\label{t-zt32}Let $f,f_{n}\in\Lambda(E)$ be such that $f(x)=\sum_{n\geq1}%
f_{n}(x)$ for every $x\in E$. Assume that the series $\sum_{n\geq1}f_{n}$
converges uniformly on a neighborhood of $x_{0}\in\operatorname*{int}%
(\operatorname*{dom}f)$. Then for every $x\in\operatorname*{int}%
(\operatorname*{dom}f)$ there exists a neighborhood $U$ of $0\in E$ with
$x+U\subset\operatorname*{dom}f$ such that the series $\sum_{n\geq1}%
f_{n+}^{\prime}(\cdot,\cdot)$ converges uniformly [to $f_{+}^{\prime}%
(\cdot,\cdot)$] on $(x+U)\times U.$
\end{theorem}

Proof. Taking into account Proposition \ref{p-zl32}, it is sufficient to prove
the conclusion for $x=x_{0}$. Moreover, as in the proof of Proposition
\ref{p-zl32}, we may (and do) assume that $x_{0}=0$ and $f_{n}(0)=0=f(0)$ for
every $n\geq1$. By hypothesis, there exists a convex neighborhood $V$ of
$x_{0}=0$ such that $2V\subset\operatorname*{dom}f$ and the series
$\sum_{n\geq1}f_{n}$ converges uniformly on $2V$. For $(x,h)\in V\times V$ we
have that $x\pm h\in2V$. It follows that the series $\sum_{n\geq1}\left[
f_{n}(x)-f_{n}(x-u)\right]  $ and $\sum_{n\geq1}\left[  f_{n}(x+u)-f_{n}%
(x)\right]  $ converge to $f(x)-f(x-u)$ and $f(x+u)-f(x)$, uniformly for
$(x,u)\in V\times V$, respectively. We have that
\[
f_{n}(x)-f_{n}(x-u)\leq\frac{f_{n}(x+tu)-f_{n}(x)}{t}\leq f_{n}(x+u)-f_{n}%
(x)\quad\forall t\in(0,1],\ \forall n\geq1.
\]
Hence
\[
\bigg\vert\sum_{k=l}^{m}\frac{f_{n}(x+tu)-f_{n}(x)}{t}\bigg\vert\leq
\bigg\vert\sum_{k=l}^{m}\left[  f_{n}(x+u)-f_{n}(x)\right]
\bigg\vert+\bigg\vert\sum_{k=l}^{m}\left[  f_{n}(x)-f_{n}(x-u)\right]
\bigg\vert
\]
for all $(x,u)\in V\times V$, all $t\in(0,1]$ and all $l,m\geq1$ with $l\leq
m$. Using (in both senses) the Cauchy criterion, we obtain that the series
$\sum_{n\geq1}\frac{f_{n}(x+tu)-f_{n}(x)}{t}$ converges uniformly for
$(x,u,t)\in V\times V\times(0,1]$ to $\frac{f(x+tu)-f(x)}{t}$. Letting
$t\rightarrow0+$, we obtain that the series $\sum_{n\geq1}f_{n+}^{\prime
}(x,u)$ converges uniformly for $(x,u)\in V\times V$ to $f_{+}^{\prime}(x,u)$.
The proof is complete. \hfill$\square$

\begin{remark}
\label{rem2}Note that for $f$, $f_{n}$ as in the preceding theorem, the series
$\sum_{n\geq1}f_{n+}^{\prime}(\cdot,\cdot)$ converges uniformly [to
$f_{+}^{\prime}(\cdot,\cdot)$] on $(x+U)\times U$ if and only if for every
$\alpha>0$, the series $\sum_{n\geq1}f_{n+}^{\prime}(\cdot,\cdot)$ converges
uniformly [to $f_{+}^{\prime}(\cdot,\cdot)$] on $(x+U)\times(\alpha U).$
\end{remark}

\begin{lemma}
\label{l1}Let $I\subset\mathbb{R}$ be an open interval and $\varphi
,\varphi_{n}:I\rightarrow\mathbb{R}$ $(n\geq1)$ be nondecreasing functions
such that $\varphi(t)=\sum_{n\geq1}\varphi_{n}(t)$ for every $t\in I$. Then
$\sum_{n\geq1}\int_{\alpha}^{\beta}\varphi_{n}(t)dt=\int_{\alpha}^{\beta
}\varphi(t)dt$ for all $\alpha,\beta\in I$ with $\alpha<\beta.$
\end{lemma}

Proof. Fix $\alpha,\beta\in I$ with $\alpha<\beta$. Take $\psi,\psi
_{n}:I\rightarrow\mathbb{R}$ $(n\geq1)$ defined by $\psi(t):=\varphi
(t)-\varphi(\alpha)$ and $\psi_{n}(t):=\varphi_{n}(t)-\varphi_{n}(\alpha)$ for
$t\in J$. Then, clearly, $\psi,\psi_{n}$ are nondecreasing functions,
$\psi_{n}(t)\geq\psi_{n}(\alpha)=0$ for $t\in J_{0}:=[\alpha,\beta]$ and
$n\geq1$ and
\[
\sum_{n\geq1}\psi_{n}(t)=\sum_{n\geq1}\left[  \varphi_{n}(t)-\varphi
_{n}(\alpha)\right]  =\sum_{n\geq1}\varphi_{n}(t)-\sum_{n\geq1}\varphi
_{n}(\alpha)=\varphi(t)-\varphi(\alpha)=\psi(t)
\]
for all $t\in I\supset J_{0}$. Since $0\leq\psi_{k}$ on $J_{0}$,
$\lim_{n\rightarrow\infty}\sum_{k=1}^{n}\psi_{k}(t)=\psi(t)$ for every $t\in
J_{0}$, and $\psi$ is Lebesgue integrable on $J_{0}$, we have that
\[
\int_{J_{0}}\psi(t)dt=\int_{J_{0}}\bigg(\sum_{n\geq1}\psi_{n}(t)\bigg)dt=\sum
_{n\geq1}\int_{J_{0}}\psi_{n}(t)dt.
\]
But $\int_{J_{0}}\psi_{n}(t)dt=\int_{J_{0}}\left(  \varphi_{n}(t)-\varphi
_{n}(\alpha)\right)  dt=\int_{\alpha}^{\beta}\varphi_{n}(t)dt-(\beta
-a)\varphi_{n}(\alpha)$, and similarly for $\psi$ and $\varphi$, whence
\begin{align*}
\int_{\alpha}^{\beta}\varphi(t)dt-(\beta-a)\varphi(\alpha)  &  =\sum_{n\geq
1}\left(  \int_{\alpha}^{\beta}\varphi_{n}(t)dt-(\beta-a)\varphi_{n}%
(\alpha)\right) \\
&  =\sum_{n\geq1}\int_{\alpha}^{\beta}\varphi_{n}(t)dt-(\beta-a)\sum_{n\geq
1}\varphi_{n}(\alpha),
\end{align*}
and so $\int_{\alpha}^{\beta}\varphi(t)dt=\sum_{n\geq1}\int_{\alpha}^{\beta
}\varphi_{n}(t)dt$. \hfill$\square$

\begin{proposition}
\label{p-zp31}Let $f,f_{n}\in\Lambda(E)$. Suppose that $\operatorname*{dom}%
f\subset\cap_{n\geq1}\operatorname*{dom}f_{n}$, and $x_{0}\in
\operatorname*{int}(\operatorname*{dom}f)$ is such that $f(x_{0})=\sum
_{n\geq1}f_{n}(x_{0}).$

\emph{(i)} Assume that $\sum_{n\geq1}f_{n+}^{\prime}(x,u)=f_{+}^{\prime}(x,u)$
for all $x\in\operatorname*{int}(\operatorname*{dom}f)$ and $u\in E$. Then
$f(x)=\sum_{n\geq1}f_{n}(x)$ for every $x\in\operatorname*{int}%
(\operatorname*{dom}f).$

\emph{(ii)} Moreover, assume that there exists a neighborhood $U$ of $0\in E$
such that $x_{0}+U\subset\operatorname*{int}(\operatorname*{dom}f)$ and the
series $\sum_{n\geq1}f_{n+}^{\prime}(\cdot,\cdot)$ converges uniformly on
$(x_{0}+U)\times U$. Then for every $x\in\operatorname*{int}%
(\operatorname*{dom}f)$ the series $\sum_{n\geq1}f_{n}$ converges uniformly to
$f$ on some neighborhood of $x.$
\end{proposition}

Proof. (i) Replacing, if necessary, $f$ by $g$ defined by $g(x):=f(x_{0}%
+x)-f(x_{0})$, and similarly for $f_{n}$, we may (an do) assume that $x_{0}=0$
and $f_{n}(x_{0})=f(0)=0$ for every $n\geq1.$

Fix $x\in\operatorname*{int}(\operatorname*{dom}f)$. Consider $I:=\{t\in
\mathbb{R}\mid tx\in\operatorname*{int}(\operatorname*{dom}f)\}$. Then $I$ is
an open interval and $0,1\in I$. Take $\theta(t):=f(tx)$ and $\theta
_{n}(t):=f_{n}(tx)$ for $t\in I$. Then $\theta,\theta_{n}$ are finite and
convex on $I$ and $\theta(0)=\theta_{n}(0)=0$. Moreover, $\theta_{+}^{\prime
}(t)=f_{+}^{\prime}(tx,x)$ for every $t\in I,$ and similarly for $\theta
_{n+}^{\prime}(t)$. Of course, $\theta$ and $\theta_{n}$ are nondecreasing and
finite on $I$. Using our hypothesis, we have that $\sum_{n\geq1}\theta
_{n+}^{\prime}(t)=\theta_{+}^{\prime}(t)$ for every $t\in I$. Using Lemma
\ref{l1} with $\varphi_{n}=\theta_{n+}^{\prime}$ and $\varphi=\theta
_{+}^{\prime}$ we obtain that
\[
f(x)=\theta(1)=\int_{0}^{1}\theta_{+}^{\prime}(t)dt=\sum_{n\geq1}\int_{0}%
^{1}\theta_{n+}^{\prime}(t)dt=\sum_{n\geq1}\theta_{n}(1)=\sum_{n\geq1}%
f_{n}(x).
\]

(ii) From (i) we have that $f(x)=\sum_{n\geq1}f_{n}(x)$ for every
$x\in\operatorname*{int}(\operatorname*{dom}f)$. Taking into account
Proposition \ref{p-zl32}, it is sufficient to show the conclusion for $x_{0}$.
As above, we may (and do) assume that $x_{0}=0$ and $f_{n}(0)=f(0)=0$. By our
hypothesis, for every $\varepsilon>0$ there exists $n_{\varepsilon}\geq0$ such
that
\[
f_{+}^{\prime}(x,u)-\varepsilon\leq\sum_{k=1}^{n}f_{k+}^{\prime}(x,u)\leq
f_{+}^{\prime}(x,u)+\varepsilon\quad\forall(x,u)\in U\times U,\ \forall n\geq
n_{\varepsilon}.
\]
In particular,
\[
f_{+}^{\prime}(tx,x)-\varepsilon\leq\sum_{k=1}^{n}f_{k+}^{\prime}(tx,x)\leq
f_{+}^{\prime}(tx,x)+\varepsilon\quad\forall t\in\lbrack0,1],\ \forall x\in
U,\ \forall n\geq n_{\varepsilon}.
\]
Integrating on $[0,1]$ with respect to $t$, we get $f(x)-\varepsilon\leq
\sum_{k=1}^{n}f_{k}(x)\leq f(x)-\varepsilon$ for all $x\in U$ and $n\geq
n_{\varepsilon}$. This shows that the conclusion holds for $x_{0}$. The proof
is complete. \hfill$\square$

\begin{theorem}
\label{t-zt33}Let $f,f_{n}\in\Lambda(E)$. Assume that $f(x)=\sum_{n\geq1}%
f_{n}(x)$ for every $x\in E$. If $f$ and $f_{n}$ are continuous on
$\operatorname*{int}(\operatorname*{dom}f)$, then
\[
\partial f(x)=w^{\ast}\text{-}\sum_{n\geq1}\partial f_{n}(x)\quad\forall
x\in\operatorname*{int}(\operatorname*{dom}f).
\]

\end{theorem}

Proof. (I) Fix $x\in\operatorname*{int}(\operatorname*{dom}f)$. Consider
$(x_{n}^{\ast})_{n\geq1}\subset\left(  \partial f_{n}(x)\right)  _{n\geq1}$.
We claim that the series $\sum_{n\geq1}x_{n}^{\ast}$ is $w^{\ast}$-convergent
to some $x^{\ast}\in\partial f(x).$

Fix some $u\in E$. By Theorem \ref{t-zt31}, for every $\varepsilon>0$, there
exists $n_{\varepsilon}\geq1$ such that $\left\vert \sum_{k=l}^{m}%
f_{k+}^{\prime}(x,\pm u)\right\vert \leq\varepsilon/2$ for all $l,m\geq
n_{\varepsilon}$ with $l\leq m$. Since $x_{k}^{\ast}\in\partial f_{k}(x)$, we
have that
\begin{equation}
\left\langle \pm u,x_{k}^{\ast}\right\rangle \leq f_{k+}^{\prime}(x,\pm u),
\label{r-z3}%
\end{equation}
and so
\[
\pm\sum_{k=l}^{m}\left\langle u,x_{k}^{\ast}\right\rangle \leq\sum_{k=l}%
^{m}f_{k+}^{\prime}(x,\pm u)\leq\bigg\vert
\sum_{k=l}^{m}f_{k+}^{\prime}(x,u)\bigg\vert +\bigg\vert
\sum_{k=l}^{m}f_{k+}^{\prime}(x,-u)\bigg\vert .
\]
Hence
\begin{equation}
\bigg\vert \sum_{k=l}^{m}\left\langle u,x_{k}^{\ast}\right\rangle
\bigg\vert \leq\bigg\vert \sum_{k=l}^{m}f_{k+}^{\prime}%
(x,u)\bigg\vert +\bigg\vert \sum_{k=l}^{m}f_{k+}^{\prime}(x,-u)\bigg\vert \leq
\varepsilon\quad\forall l,m\geq1,\ n_{\varepsilon}\leq l\leq m. \label{r-z4}%
\end{equation}
Therefore, the series $\sum_{n\geq1}\left\langle u,x_{n}^{\ast}\right\rangle $
is convergent, and so $\varphi(u):=\sum_{n\geq1}\left\langle u,x_{k}^{\ast
}\right\rangle \in\mathbb{R}$. Moreover, from (\ref{r-z3}) and Theorem
\ref{t-zt31} we get
\[
\varphi(u)\leq\sum_{n\geq1}f_{n+}^{\prime}(x,u)=f_{+}^{\prime}(x,u).
\]
We got so a linear mapping $\varphi:E\rightarrow\mathbb{R}$ such that
$\varphi\leq f_{+}^{\prime}(x,\cdot)$. Since $f$ is continuous at
$x\in\operatorname*{dom}f$, $f_{+}^{\prime}(x,\cdot)$ is continuous, and so
$\varphi\in\partial f(x)$. Hence $w^{\ast}$-$\sum_{n\geq1}x_{n}^{\ast}$ exists
and belongs to $\partial f(x)$. Therefore, condition (I) in Definition
\ref{d-z1} holds.

(II) Taking the limit for $m\rightarrow\infty$ in (\ref{r-z4}) we obtain that
$\big\vert \sum_{k\geq n}\left\langle u,x_{k}^{\ast}\right\rangle
\big\vert \leq\varepsilon$ for all $n\geq n_{\varepsilon}$. Since
$n_{\varepsilon}$ does not depend on the sequence $(x_{n}^{\ast})_{n\geq
1}\subset\left(  \partial f_{n}(x)\right)  _{n\geq1}$, the second condition in
Definition \ref{d-z1} holds, too.

(III) For $n\geq1$ set $R_{n}:=\sum_{k\geq n}f_{k}$; clearly $R_{n}\in
\Lambda(E)$ for $n\geq1$. Using Theorem \ref{t-zt31} for the sequence
$(f_{k+n})_{k\geq0}$ we obtain that $R_{n+}^{\prime}(x,u)=\sum_{k\geq n}%
f_{k+}^{\prime}(x,u)\rightarrow0$ (as $n\rightarrow\infty$) for every $u\in
E$. Moreover, $R_{1}=f$ and $R_{k}=f_{k}+R_{k+1}$ for $k\geq1$. Since $f_{k}$
is continuous on $\operatorname*{int}(\operatorname*{dom}f_{k})\supset
\operatorname*{int}(\operatorname*{dom}f)=\operatorname*{int}%
(\operatorname*{dom}R_{k})$, by Moreau-Rockafellar Theorem one has
\begin{equation}
\partial R_{k}(x)=\partial f_{k}(x)+\partial R_{k+1}(x)\quad\forall k\geq1.
\label{r-z5}%
\end{equation}
Consider $x^{\ast}\in\partial f(x)=\partial R_{1}(x)$. From (\ref{r-z5})
applied for $k=1$ we get $x_{1}^{\ast}\in\partial f_{1}(x)$ such that
$x^{\ast}-x_{1}^{\ast}\in\partial R_{2}(x).$ Continuing in this way we get a
sequence $(x_{n}^{\ast})_{n\geq1}\subset(\partial f_{n}(x))_{n\geq1}$ such
that
\[
y_{n}^{\ast}:=x^{\ast}-\sum_{k=1}^{n-1}x_{k}^{\ast}\in\partial R_{n}%
(x)\quad\forall n\geq1,
\]
where $\sum_{k=1}^{0}x_{k}^{\ast}:=0$. Using Theorem \ref{t-zt31} for the
sequence $(f_{k+n})_{k\geq0}$ and the fact that $\left\langle \pm
u,y_{n}^{\ast}\right\rangle \leq R_{n+}^{\prime}(x,\pm u)$ we obtain that
\[
\bigg\vert\bigg\langle u,x^{\ast}-\sum_{k=1}^{n}x_{k}^{\ast}%
\bigg\rangle\bigg\vert=\left\vert \left\langle u,y_{n}^{\ast}\right\rangle
\right\vert \leq\left\vert R_{n+}^{\prime}(x,u)\right\vert +\left\vert
R_{n+}^{\prime}(x,-u)\right\vert \rightarrow0~\text{for }n\rightarrow\infty
\]
for every $u\in E$. It follows that $x^{\ast}=w^{\ast}$-$\sum_{n\geq1}%
x_{n}^{\ast}$. Hence condition (III) of Definition \ref{d-z1} is verified,
too. The proof is complete. \hfill$\square$

\medskip

When $E$ is a normed vector space, one has also the next result.

\begin{theorem}
\label{t-zt34}Let $E$ be a normed vector space and $f,f_{n}\in\Lambda(E).$
Assume that $f(x)=\sum_{n\geq1}f_{n}(x)$ for every $x\in E$. If $f$ and
$f_{n}$ are continuous on $\operatorname*{int}(\operatorname*{dom}f)$ and the
series $\sum_{n\geq1}f_{n}$ converges uniformly on a nonempty open subset of
$\operatorname*{dom}f$, then
\[
\partial f(x)=\left\Vert \cdot\right\Vert \text{-}\sum_{n\geq1}\partial
f_{n}(x)\quad\forall x\in\operatorname*{int}(\operatorname*{dom}f);
\]
moreover, $\lim_{n\rightarrow\infty}\big\Vert \sum_{k\geq n}\partial
f_{k}(x)\big\Vert =0$ uniformly on some neighborhood of $x$ for every
$x\in\operatorname*{int}(\operatorname*{dom}f)$, where $\left\Vert
A\right\Vert :=\sup\left\{  \left\Vert x^{\ast}\right\Vert \mid x^{\ast}\in
A\right\}  $ for $\emptyset\neq A\subset X^{\ast}.$
\end{theorem}

Proof. We follow the same steps as in the proof of Theorem \ref{t-zt33}. So,
fix $x_{0}\in\operatorname*{int}(\operatorname*{dom}f).$

(I) Consider $(x_{n}^{\ast})_{n\geq1}\subset\left(  \partial f_{n}%
(x_{0})\right)  _{n\geq1}$. By Theorem \ref{t-zt33} we have that $w^{\ast}%
$-$\sum_{n\geq1}x_{n}^{\ast}=x^{\ast}\in\partial f(x_{0})$. We claim that
$x^{\ast}=\left\Vert \cdot\right\Vert $-$\sum_{n\geq1}x_{n}^{\ast}$. Indeed,
using this time Proposition \ref{p-zl32}, Theorem \ref{t-zt32} and Remark
\ref{rem2}, there exists $r>0$ such that the series $\sum_{n\geq1}%
f_{n+}^{\prime}(\cdot,\cdot)$ converges uniformly to $f_{+}^{\prime}%
(\cdot,\cdot)$ on $(x_{0}+rU_{E})\times U_{E}$, where $U_{E}$ is the closed
unit ball of $E$. Taking $\varepsilon>0$, as in the proof of Theorem
\ref{t-zt33}, there exists $n_{\varepsilon}\geq1$ such that (\ref{r-z4}) holds
for all $(x,u)\in(x_{0}+rU_{E})\times U_{E}$, and so $\left\Vert \sum
_{k=l}^{m}x_{k}^{\ast}\right\Vert \leq\varepsilon$ for all $n_{\varepsilon
}\leq l\leq m$. Since $x^{\ast}=w^{\ast}$-$\sum_{n\geq1}x_{n}^{\ast}$ we get
$x^{\ast}=\left\Vert \cdot\right\Vert $-$\sum_{n\geq1}x_{n}^{\ast}.$ In fact
we even get $\lim_{n\rightarrow\infty}\big\Vert \sum_{k\geq n}\partial
f_{k}(x)\big\Vert =0$ uniformly on $x_{0}+rU_{E}.$

(II) This step is practically stated in (I).

(III) Having $x^{\ast}\in\partial f(x_{0})$, from Theorem \ref{t-zt33} we find
$(x_{n}^{\ast})_{n\geq1}\subset\left(  \partial f_{n}(x_{0})\right)  _{n\geq
1}$ such that $x^{\ast}=w^{\ast}$-$\sum_{n\geq1}x_{n}^{\ast}$. From (I) we
obtain that $x^{\ast}=\left\Vert \cdot\right\Vert $-$\sum_{n\geq1}x_{n}^{\ast
}$. The proof is complete. \hfill$\square$

\begin{corollary}
\label{c-zc31}Let $f,f_{n}\in\Lambda(E)$. Assume that $f(x)=\sum_{n\geq1}%
f_{n}(x)$ for every $x\in E$, and $f$, $f_{n}$ are continuous on
$\operatorname*{int}(\operatorname*{dom}f)$ for every $n\geq1$. Take
$\overline{x}\in\operatorname*{int}(\operatorname*{dom}f)$. Then

\emph{(i)} $f$ is G\^{a}teaux differentiable at $\overline{x}$ if and only if
$f_{n}$ is G\^{a}teaux differentiable at $\overline{x}$ for every $n\geq1$, in
which case $\nabla f(\overline{x})=\sum_{n\geq1}\nabla f_{n}(\overline{x}).$

\emph{(ii)} Moreover, assume that $E$ is a normed vector space. If $f$ is
Fr\'{e}chet differentiable at $\overline{x}$ then $f_{n}$ is Fr\'{e}chet
differentiable at $\overline{x}$ for every $n\geq1.$
\end{corollary}

Proof. By Theorem \ref{t-zt33} we have that $\partial f(\overline{x})=w^{\ast
}$-$\sum_{n\geq1}\partial f_{n}(\overline{x})$. This relation shows that
$\partial f(\overline{x})$ is a singleton if and only if $\partial
f_{n}(\overline{x})$ is a singleton for every $n\geq1$. Since the functions
$f$ and $f_{n}$ $(n\geq1)$ are continuous at $\overline{x}$, (i) follows.

(ii) Assume now that $E$ is a nvs and $f$ is Fr\'{e}chet differentiable at
$\overline{x}$. It is known that $g$ and $h$ are Fr\'{e}chet differentiable at
$\overline{x}\in\operatorname*{int}(\operatorname*{dom}g\cap
\operatorname*{dom}h)=\operatorname*{int}(\operatorname*{dom}(g+h))$ provided
$g,h\in\Lambda(E)$, $g,h$ are continuous at $\overline{x}$ and $g+h$ is
Fr\'{e}chet differentiable at $\overline{x}$. This is due to the fact that in
such conditions, as seen from (i), $g$ and $h$ are G\^{a}teaux differentiable
at $\overline{x}$. Then we have%
\[
0\leq\frac{g(\overline{x}+u)-g(\overline{x})-\left\langle u,\nabla
g(\overline{x})\right\rangle }{\left\Vert u\right\Vert }\leq\frac
{(g+h)(\overline{x}+u)-(g+h)(\overline{x})-\left\langle u,\nabla
(g+h)(\overline{x})\right\rangle }{\left\Vert u\right\Vert }\rightarrow0
\]
for $\left\Vert u\right\Vert \rightarrow0.$

With the notation in the proof of Theorem \ref{t-zt33}, $f=R_{1}=f_{1}+R_{2}.$
It follows that $f_{1}$ and $R_{2}$ are Fr\'{e}chet differentiable at
$\overline{x}$. Since $R_{2}=f_{2}+R_{3},$ it follows that $f_{2}$ and $R_{3}$
are Fr\'{e}chet differentiable at $\overline{x}$. Continuing in this way we
obtain that $f_{n}$ is Fr\'{e}chet differentiable at $\overline{x}$ for every
$n\geq1.$ \hfill$\square$

\medskip

Of course, if $\dim E<\infty$, the weak$^{\ast}$ and strong convergences on
$E^{\ast}$ coincide, and so Theorems \ref{t-zt33} and \ref{t-zt34} are
equivalent; moreover, G\^{a}teaux and Fr\'{e}chet differentiability for convex
functions coincide, and so the converse implication in Corollary \ref{c-zc31}
(ii) is true.

\begin{question}
\textbf{\label{q1}}Is the converse of Corollary \ref{c-zc31} (ii) true when
$\dim E=\infty$?
\end{question}

In the sequel we set $\mathbb{R}_{+}:=[0,\infty)$, $\mathbb{R}_{+}^{\ast
}:=(0,\infty),$ $\mathbb{R}_{-}:=-\mathbb{R}_{+}$, $\mathbb{R}_{-}^{\ast
}:=-\mathbb{R}_{+}^{\ast}.$

\begin{proposition}
\label{prop1}Let $f_{n}(x):=e^{\sigma_{n}x}$ for $x\in\mathbb{R}$ with
$(\sigma_{n})_{n\geq1}\subset\mathbb{R}$; set $f=\sum_{n\geq1}f_{n}.$

\emph{(i)} If $\overline{x}\in\operatorname*{dom}f$ then $\sigma_{n}%
\overline{x}\rightarrow-\infty$, and so either $\overline{x}>0$ and
$\sigma_{n}\rightarrow-\infty$, or $\overline{x}<0$ and $\sigma_{n}%
\rightarrow\infty.$

\emph{(ii)} Assume that $0<\sigma_{n}\rightarrow\infty$ and
$\operatorname*{dom}f\neq\emptyset$. Then there exists $\alpha\in
\mathbb{R}_{+}$ such that $I:=(-\infty,-\alpha)\subset\operatorname*{dom}%
f\subset\operatorname*{cl}I$, $f$ is strictly convex and increasing on
$\operatorname*{dom}f$, and $\lim_{x\rightarrow-\infty}f(x)=0=\inf f.$
Moreover,
\begin{equation}
f^{\prime}(x)=\sum_{n\geq1}f_{n}^{\prime}(x)=\sum_{n\geq1}\sigma_{n}%
e^{\sigma_{n}x}\quad\forall x\in\operatorname*{int}(\operatorname*{dom}f)=I,
\label{r-8}%
\end{equation}
$f^{\prime}$ is increasing and continuous on $I$, $\lim_{x\rightarrow-\infty
}f^{\prime}(x)=0$, and
\begin{equation}
\lim_{x\uparrow-\alpha}f^{\prime}(x)=\sum_{n\geq1}\sigma_{n}e^{-\sigma
_{n}\alpha}=:\gamma\in(0,\infty]. \label{r-7}%
\end{equation}
In particular, $\partial f(\operatorname*{int}(\operatorname*{dom}%
f))=f^{\prime}(I)=(0,\gamma).$

\emph{(iii)} Assume that $0<\sigma_{n}\rightarrow\infty$ is such that
$\operatorname*{int}(\operatorname*{dom}f)=I:=(-\infty,-\alpha)$ with
$\alpha\in\mathbb{R}_{+}^{\ast}$. Then either \emph{(a) }$\operatorname*{dom}%
f=I$ and $\gamma=\infty$, or \emph{(b)} $\operatorname*{dom}%
f=\operatorname*{cl}I$ and $\gamma=\infty$, in which case $f_{-}^{\prime
}(-\alpha)=\gamma$, $\partial f(-\alpha)=\emptyset$ and the series
$\sum_{n\geq1}f_{n}^{\prime}(-\alpha)$ is not convergent, or \emph{(c)}
$\operatorname*{dom}f=\operatorname*{cl}I$ and $\gamma<\infty$, in which case
$f_{-}^{\prime}(-\alpha)=\gamma$ and%
\[
\partial f(-\alpha)=[\gamma,\infty)\neq\left\{  \gamma\right\}  =\sum_{n\geq
1}\partial f_{n}(-\alpha).
\]

\end{proposition}

Proof. (i) Take $\overline{x}\in\operatorname*{dom}f$. Then the series
$\sum_{n\geq1}e^{\sigma_{n}\overline{x}}$ is convergent, and so $e^{\sigma
_{n}\overline{x}}\rightarrow0$. The conclusion is obvious.

(ii) Set $\beta:=\sup(\operatorname*{dom}f)\in(-\infty,\infty]$ and take
$x<\beta$. Then there exists $\overline{x}\in\operatorname*{dom}f$ with
$x\leq\overline{x}.$ Because $\sigma_{n}\geq0$, and so $0<e^{\sigma_{n}x}\leq
e^{\sigma_{n}\overline{x}}$ for $n\geq1$, the series $\sum_{n\geq1}%
e^{\sigma_{n}x}$ is convergent, whence $x\in\operatorname*{dom}f$. Hence
$(-\infty,\beta)\subset\operatorname*{dom}f\subset(-\infty,\beta].$ Since
$0\notin\operatorname*{dom}f$, we obtain that $\beta\leq0,$ and so
$\alpha:=-\beta$ does the job. Since $f_{n}$ is strictly convex and
increasing, $f$ is strictly convex and increasing on its domain.

Because $f_{n}^{\prime}(x)=\sigma_{n}e^{\sigma_{n}x}$ for every $x\in
\mathbb{R}$, we get (\ref{r-8}) using Corollary \ref{c-zc31}. From (\ref{r-8})
we have that $f^{\prime}$ is increasing and continuous on $\operatorname*{int}%
(\operatorname*{dom}f).$

Since $0<\sigma_{n}e^{\sigma_{n}x}\leq\sigma_{n}e^{\sigma_{n}\overline{x}}$
for all $n\geq1$ and $x\leq\overline{x}$, the series $\sum_{n\geq1}\sigma
_{n}e^{\sigma_{n}x}$ is uniformly convergent (u.c.~for short) on
$(-\infty,\overline{x}]$ for any $\overline{x}\in\operatorname*{int}%
(\operatorname*{dom}f)$. Since $\lim_{x\rightarrow-\infty}\left(  \sigma
_{n}e^{\sigma_{n}x}\right)  =\lim_{x\rightarrow-\infty}e^{\sigma_{n}x}=0$,
from (\ref{r-8}) and $f=\sum_{n\geq1}f_{n}$ we obtain that $\lim
_{x\rightarrow-\infty}f^{\prime}(x)=0$ and $\lim_{x\rightarrow-\infty}f(x)=0$, respectively.

From (\ref{r-8}) we have that
\begin{equation}
f^{\prime}(x)=\sum_{n\geq1}\sigma_{n}e^{\sigma_{n}x}\geq\sum_{k=1}^{n}%
\sigma_{k}e^{\sigma_{k}x}\quad\forall x\in I,~\forall n\geq1. \label{r-9}%
\end{equation}
Taking the limit for $-\alpha>x\rightarrow-\alpha$, we get $\lim
_{x\uparrow-\alpha}f^{\prime}(x)\geq\sum_{k=1}^{n}\sigma_{k}e^{-\sigma
_{k}\alpha}$. Taking now the limit for $n\rightarrow\infty$ we get
$\lim_{x\uparrow-\alpha}f^{\prime}(x)\geq\gamma:=\sum_{n\geq1}\sigma
_{n}e^{-\sigma_{n}\alpha}\in(0,\infty]$. If $\gamma=\infty$ it is clear that
(\ref{r-7}) holds. Assume that $\gamma<\infty$. There exists some $n_{0}\geq1$
such that $\sigma_{n}\geq1$, whence $\sigma_{n}e^{-\sigma_{n}\alpha}\geq
e^{-\sigma_{n}\alpha}$, for $n\geq n_{0}$, and so $\sum_{n\geq1}e^{-\sigma
_{n}\alpha}$ is convergent. It follows that the series $\sum_{n\geq1}f_{n}$
and $\sum_{n\geq1}f_{n}^{\prime}$ are u.c.~on $(-\infty,-\alpha]$, and so
$\lim_{x\uparrow-\alpha}f^{\prime}(x)=\sum_{n\geq1}\lim_{x\uparrow-\alpha
}f_{n}^{\prime}(x)=\sum_{n\geq1}\sigma_{n}e^{-\sigma_{n}\alpha},$ that is
(\ref{r-7}) holds in this case, too.

(iii) Let $\alpha\in\mathbb{R}_{+}^{\ast}$. Since $0<e^{-\sigma_{n}\alpha}%
\leq\sigma_{n}e^{-\sigma_{n}\alpha}$ for large $n$, we get $\gamma=\infty$
when $-\alpha\notin\operatorname*{dom}f$, and so (a) holds. Assume that
$-\alpha\in\operatorname*{dom}f$. Since $f\in\Gamma(\mathbb{R})$, we have that
$f_{-}^{\prime}$ is continuous from the left, whence $f_{-}^{\prime}%
(-\alpha)=\lim_{x\uparrow-\alpha}f^{\prime}(x)$, and $\partial f(-\alpha
)=[f_{-}^{\prime}(-\alpha),\infty)$. Now the conclusion is immediate using
(\ref{r-7}). \hfill$\square$

\begin{example}
\label{ex1}In Proposition \ref{prop1}, for $\sigma_{n}=n^{\theta}$ $(n\geq1)$
with $\theta>0$ one has $\operatorname*{dom}f=(-\infty,0)$, for $\sigma
_{n}=\ln\left[  n(\ln n)^{\theta}\right]  $ $(n\geq2)$ with $\theta
\in\mathbb{R}$ one has $\operatorname*{int}(\operatorname*{dom}f)=(-\infty
,-1)$, while for $\sigma_{n}=\ln(\ln n)$ $(n\geq2)$ one has
$\operatorname*{dom}f=\emptyset$. Moreover, let $\sigma_{n}=\ln\left[  n(\ln
n)^{\theta}\right]  $ $(n\geq2)$;\footnote{We thank Prof.\ C. Lefter for this
example.} for $\theta\in(-\infty,1]$ one has $\operatorname*{dom}%
f=(-\infty,-1)$, for $\theta\in(1,2]$ one has $\operatorname*{dom}%
f=(-\infty,-1]$ and $f_{-}^{\prime}(-1)=\infty$, for $\theta\in(2,\infty)$ one
has $\operatorname*{dom}f=(-\infty,-1]$ and $f_{-}^{\prime}(-1)<\infty.$
\end{example}

Proposition \ref{prop1} (iii) (b) and Example \ref{ex1} show that the
conclusion of Theorem \ref{t-zt33} can be false for $x\in\operatorname*{dom}%
f\setminus\operatorname*{int}(\operatorname*{dom}f)$ [even for $x\in
\operatorname*{dom}(\partial f)\setminus\operatorname*{int}%
(\operatorname*{dom}f)$].

One could ask if the condition $\operatorname*{int}(\operatorname*{dom}%
f)\neq\emptyset$ in Theorem \ref{t-zt33} is just a technical assumption. The
next example shows that this condition is essential.

\begin{example}
\label{ex2}Let $g_{n}(x,y):=e^{nx+(-1)^{n}\varsigma_{n}y}$ for $x,y\in
\mathbb{R}$ and $g=\sum_{n\geq1}g_{n}$, where $(\varsigma_{n})_{n\geq1}%
\subset\mathbb{R}$ is such that $\varsigma_{n}/n\rightarrow\infty$. Clearly
$g,g_{n}\in\Gamma(\mathbb{R}^{2})$ with $g(x,y)=f(x)+\iota_{\{0\}}(y)$ for
$(x,y)\in\mathbb{R}^{2}$, where
\begin{equation}
f(x)=\sum_{n\geq1}e^{nx}=\left\{
\begin{array}
[c]{ll}%
\frac{e^{x}}{1-e^{x}} & \text{if }x\in\mathbb{R}_{-}^{\ast},\\
\infty & \text{if }x\in\mathbb{R}_{+}.
\end{array}
\right.  \label{r-2}%
\end{equation}
Hence $\operatorname*{dom}g=\mathbb{R}_{-}^{\ast}\times
\{0\}=\operatorname*{ri}(\operatorname*{dom}g)$, but $\operatorname*{int}%
(\operatorname*{dom}g)=\emptyset$. It is clear that for $(x,y)\in
\operatorname*{ri}(\operatorname*{dom}g)=\mathbb{R}_{-}^{\ast}\times\{0\}$ we
have that $\partial g(x,0)=\partial f(x)\times\partial\iota_{\{0\}}%
(0)=\{e^{x}/(1-e^{x})^{2}\}\times\mathbb{R}$. However, $\partial
g_{n}(x,0)=\left\{  \nabla g_{n}(x,0)\right\}  =\{(ne^{nx},(-1)^{n}%
\varsigma_{n}e^{nx})\}$. For $\varsigma_{n}:=n^{2}$ $(n\geq1)$ we get
\begin{equation}
\sum_{n\geq1}\nabla g_{n}(x,0)=\big(f^{\prime}(x),8f^{\prime\prime
}(2x)-f^{\prime\prime}(x)\big)\quad\forall x\in\mathbb{R}_{-}^{\ast},
\label{r5}%
\end{equation}
for $\varsigma_{n}=e^{\alpha n}$ with $\alpha>0$ we get $\sum_{n\geq1}\nabla
g_{n}(x,0)=\big(f^{\prime}(x),-e^{x+\alpha}/\left(  1+e^{x+\alpha}\right)
\big)$ for $x<-\alpha$ and $\sum_{n\geq1}\nabla g_{n}(x,0)$ is not convergent
for $x\in\lbrack-\alpha,0)$, while for $\varsigma_{n}=e^{n^{2}}$ $(n\geq1)$
the series $\sum_{n\geq1}\nabla g_{n}(x,0)$ is not convergent for each
$x\in\mathbb{R}_{-}^{\ast}$.

Indeed, we have that $f^{\prime}(x)=\sum_{n\geq1}ne^{nx}=e^{x}/(1-e^{x})^{2}$
and $f^{\prime\prime}(x)=\sum_{n\geq1}n^{2}e^{nx}$ for $x\in\mathbb{R}%
_{-}^{\ast}$. It follows that for $x\in\mathbb{R}_{-}^{\ast}$ we have
\[
8f^{\prime\prime}(2x)=2\sum_{n\geq1}(2n)^{2}e^{2nx}=\sum_{n\geq1}(-1)^{n}%
n^{2}e^{nx}+\sum_{n\geq1}n^{2}e^{nx}=\sum_{n\geq1}(-1)^{n}n^{2}e^{nx}%
+f^{\prime\prime}(x),
\]
whence (\ref{r5}) follows for $\varsigma_{n}=n^{2}.$
\end{example}

\section{Applications to the conjugate of a countable sum}

A natural question is what we could say about the conjugate of $f=\sum
_{n\geq1}f_{n}$ when $f,f_{n}\in\Lambda(E)$. It is known that for a finite
family of functions $f_{1},\ldots,f_{n}\in\Lambda(E)$ one has
\begin{align*}
(f_{1}+\ldots+f_{n})^{\ast}(x^{\ast})  &  \leq\left(  f_{1}^{\ast}%
\square\ldots\square f_{n}^{\ast}\right)  (x^{\ast})\\
&  :=\inf\left\{  f_{1}^{\ast}(x_{1}^{\ast})+\ldots+f_{n}^{\ast}(x_{n}^{\ast
})\mid x_{1}^{\ast}+\ldots+x_{n}^{\ast}=x^{\ast}\right\}  .
\end{align*}
Of course, in the above formula one could take $x_{k}^{\ast}\in
\operatorname*{dom}f_{k}^{\ast}$ for every $k\in\overline{1,n}$ using the
usual convention $\inf\emptyset:=\infty$. Moreover, when all functions (but
one) are continuous at some point in $\cap_{k\in\overline{1,n}}%
\operatorname*{dom}f_{k}$ we have, even for every $x^{\ast}\in E^{\ast},$%
\[
(f_{1}+\ldots+f_{n})^{\ast}(x^{\ast})=\min\left\{  f_{1}^{\ast}(x_{1}^{\ast
})+\ldots+f_{n}^{\ast}(x_{n}^{\ast})\mid x_{1}^{\ast}+\ldots+x_{n}^{\ast
}=x^{\ast}\right\}  .
\]
Recall that the inf-convolution operation $\square$ was introduced by J. J.
Moreau in \cite{Mor:63}; many properties of this operation can be found in
\cite{Mor:66}, among them being the formula mentioned above. The aim of this
section is to extend and study this operation to countable sums of convex functions.

In the next proposition we put together several assertions on the conjugate of
$\sum_{n\geq1}f_{n}$; the last assertion is an application of Theorem 9.

\begin{proposition}
\label{p3}Let $f,f_{n}\in\Lambda(E)$. Assume that $f(x)=\sum_{n\geq1}f_{n}(x)$
for every $x\in E.$

\emph{(i)} If $(x_{n}^{\ast})_{n\geq1}\subset(\operatorname*{dom}f_{n}^{\ast
})_{n\geq1}$ is such that $w^{\ast}$-$\sum_{n\geq1}x_{n}^{\ast}=x^{\ast}\in
E^{\ast}$, then the sequence $\sum_{k=1}^{n}f_{k}^{\ast}(x_{k}^{\ast})$ has a
limit in $(-\infty,\infty]$ and $f^{\ast}(x^{\ast})\leq\sum_{n\geq1}%
f_{n}^{\ast}(x_{n}^{\ast})$; in particular,
\begin{equation}
f^{\ast}(x^{\ast})\leq\inf\bigg\{\sum_{n\geq1}f_{n}^{\ast}(x_{n}^{\ast}%
)\mid(x_{n}^{\ast})_{n\geq1}\subset(\operatorname*{dom}f_{n}^{\ast})_{n\geq
1},\ x^{\ast}=w^{\ast}\text{-}\sum_{n\geq1}x_{n}^{\ast}\bigg\} \label{r-1}%
\end{equation}
for every $x^{\ast}\in E^{\ast}$.

\emph{(ii)} If there exists $(x_{n}^{\ast})_{n\geq1}\subset
(\operatorname*{dom}f_{n}^{\ast})_{n\geq1}$ such that $w^{\ast}$-$\sum
_{n\geq1}x_{n}^{\ast}=x^{\ast}\in E^{\ast}$ and $\sum_{n\geq1}f_{n}^{\ast
}(x_{n}^{\ast})\in\mathbb{R}$ then $x^{\ast}\in\operatorname*{dom}f^{\ast}.$

\emph{(iii)} If $\overline{x}\in\cap_{n\geq1}\operatorname*{dom}f_{n}$ and
$(\overline{x}_{n}^{\ast})_{n\geq1}\subset\left(  \partial f_{n}(\overline
{x})\right)  _{n\geq1}$ is such that $w^{\ast}$-$\sum_{n\geq1}\overline{x}%
_{n}^{\ast}=x^{\ast}\in E^{\ast}$, then $\overline{x}\in\operatorname*{dom}f$,
$x^{\ast}\in\partial f(\overline{x})$ and $f^{\ast}(x^{\ast})=\sum_{n\geq
1}f_{n}^{\ast}(\overline{x}_{n}^{\ast})$. In particular,%
\begin{equation}
f^{\ast}(x^{\ast})=\min\bigg\{\sum_{n\geq1}f_{n}^{\ast}(x_{n}^{\ast}%
)\mid(x_{n}^{\ast})_{n\geq1}\subset(\operatorname*{dom}f_{n}^{\ast})_{n\geq
1},\ x^{\ast}=w^{\ast}\text{-}\sum_{n\geq1}x_{n}^{\ast}\bigg\} \label{r2}%
\end{equation}

\emph{(iv)} Let $x^{\ast}\in\partial f(\overline{x})$ with $\overline{x}%
\in\operatorname*{dom}f$ and let $(x_{n}^{\ast})_{n\geq1}\subset
(\operatorname*{dom}f_{n}^{\ast})_{n\geq1}$ be such that $x^{\ast}=w^{\ast}%
$-$\sum_{n\geq1}x_{n}^{\ast}$. Then $f^{\ast}(x^{\ast})=\sum_{n\geq1}%
f_{n}^{\ast}(x_{n}^{\ast})$ if and only if $x_{n}^{\ast}\in\partial
f_{n}(\overline{x})$ for every $n\geq1.$

\emph{(v)} Assume that $f$ and $f_{n}$ $(n\geq1)$ are continuous on
$\operatorname*{int}(\operatorname*{dom}f)$. Then for every $x^{\ast}%
\in\partial f\left(  \operatorname*{int}(\operatorname*{dom}f)\right)  $
relation (\ref{r2}) holds. More precisely, if $x^{\ast}\in\partial f(x)$ for
$x\in\operatorname*{int}(\operatorname*{dom}f)$ then $x^{\ast}=w^{\ast}$%
-$\sum_{n\geq1}x_{n}^{\ast}$ for some $(x_{n}^{\ast})_{n\geq1}\subset(\partial
f_{n}(x))_{n\geq1}$, and $f^{\ast}(x^{\ast})=\sum_{n\geq1}f_{n}^{\ast}%
(x_{n}^{\ast}).$
\end{proposition}

Proof. (i) Let us fix $(x_{n}^{\ast})_{n\geq1}\subset(\operatorname*{dom}%
f_{n}^{\ast})_{n\geq1}$ with $w^{\ast}$-$\sum_{n\geq1}x_{n}^{\ast}=x^{\ast}\in
X^{\ast}$; take $x\in\operatorname*{dom}f$ and $\gamma_{n}:=f_{n}%
(x)+f_{n}^{\ast}(x_{n}^{\ast})-\left\langle x,x_{n}^{\ast}\right\rangle \geq0$
for $n\geq1$. Hence $\sum_{k=1}^{n}\gamma_{k}\rightarrow\gamma\in
\lbrack0,\infty]$, and so
\begin{align*}
\sum_{n\geq1}f_{n}^{\ast}(x_{n}^{\ast})  &  =\lim_{n\rightarrow\infty}%
\sum_{k=1}^{n}f_{k}^{\ast}(x_{k}^{\ast})=\lim_{n\rightarrow\infty}%
\bigg(\sum_{k=1}^{n}\gamma_{k}-\sum_{k=1}^{n}f_{k}(x)+\bigg\langle x,\sum
_{k=1}^{n}x_{k}^{\ast}\bigg\rangle\bigg)\\
&  =\gamma-f(x)+\left\langle x,x^{\ast}\right\rangle \geq\left\langle
x,x^{\ast}\right\rangle -f(x).
\end{align*}
It follows that $\sum_{n\geq1}f_{n}^{\ast}(x_{n}^{\ast})\geq\sup
_{x\in\operatorname*{dom}f}\left(  \left\langle x,x^{\ast}\right\rangle
-f(x)\right)  =f^{\ast}(x^{\ast})>-\infty$. Relation (\ref{r-1}) is now obvious.

(ii) The assertion is an immediate consequence of (\ref{r-1}).

(iii) Since $\overline{x}_{n}^{\ast}\in\partial f_{n}(\overline{x}%
)\subset\operatorname*{dom}f_{n}^{\ast}$, we have that $f_{n}^{\ast}%
(\overline{x}_{n}^{\ast})=\left\langle \overline{x},\overline{x}_{n}^{\ast
}\right\rangle -f_{n}(\overline{x})$ for $n\geq1$, and so, using (i), we get%
\[
-\infty<f^{\ast}(x^{\ast})\leq\sum_{n\geq1}f_{n}^{\ast}(\overline{x}_{n}%
^{\ast})=\sum_{n\geq1}\left[  \left\langle \overline{x},\overline{x}_{n}%
^{\ast}\right\rangle -f_{n}(\overline{x})\right]  =\left\langle \overline
{x},x^{\ast}\right\rangle -f(\overline{x})<\infty.
\]
It follows that $\overline{x}\in\operatorname*{dom}f$ and $f^{\ast}(x^{\ast
})+f(\overline{x})\leq\left\langle \overline{x},x^{\ast}\right\rangle ,$
whence $x^{\ast}\in\partial f(\overline{x})$ and $f^{\ast}(x^{\ast}%
)=\sum_{n\geq1}f_{n}^{\ast}(\overline{x}_{n}^{\ast})$. Using again (i) we
obtain that (\ref{r2}) holds.

(iv) Since $\operatorname*{dom}f\subset\cap_{n\geq1}\operatorname*{dom}f_{n},$
we have that $\overline{x}\in\cap_{n\geq1}\operatorname*{dom}f_{n}$. Assuming
that $x_{n}^{\ast}\in\partial f_{n}(\overline{x})$ for $n\geq1$, the
conclusion $f^{\ast}(x^{\ast})=\sum_{n\geq1}f_{n}^{\ast}(x_{n}^{\ast})$
follows from (iii).

Conversely, assume that $f^{\ast}(x^{\ast})=\sum_{n\geq1}f_{n}^{\ast}%
(x_{n}^{\ast})$. Then
\[
0=\sum_{n\geq1}f_{n}^{\ast}(x_{n}^{\ast})+\sum_{n\geq1}f_{n}(\overline
{x})-\sum_{n\geq1}\left\langle \overline{x},x_{n}^{\ast}\right\rangle
=\sum_{n\geq1}\left[  f_{n}^{\ast}(x_{n}^{\ast})+f_{n}(\overline
{x})-\left\langle \overline{x},x_{n}^{\ast}\right\rangle \right]  .
\]
Since $f_{n}^{\ast}(x_{n}^{\ast})+f_{n}(\overline{x})-\left\langle
\overline{x},x_{n}^{\ast}\right\rangle \geq0$ for $n\geq1$, we obtain that
$f_{n}^{\ast}(x_{n}^{\ast})+f_{n}(\overline{x})-\left\langle \overline
{x},x_{n}^{\ast}\right\rangle =0$, and so $x_{n}^{\ast}\in\partial
f_{n}(\overline{x})$ for every $n\geq1.$

(v) Assume now that $f$ and $f_{n}$ are continuous on $\operatorname*{int}%
(\operatorname*{dom}f)$ and take $x^{\ast}\in\partial f\left(
\operatorname*{int}(\operatorname*{dom}f)\right)  $. Then there exists
$\overline{x}\in\operatorname*{int}(\operatorname*{dom}f)$ such that $x^{\ast
}\in\partial f(\overline{x})$. By Theorem \ref{t-zt33}, there exists
$(x_{n}^{\ast})_{n\geq1}\subset(\partial f_{n}(\overline{x}))_{n\geq1}$ such
that $x^{\ast}=w^{\ast}$-$\sum_{n\geq1}x_{n}^{\ast}$. By (iii) we have that
$f^{\ast}(x^{\ast})=\sum_{n\geq1}f_{n}^{\ast}(x_{n}^{\ast}).$ \hfill$\square$

\begin{remark}
\label{rem1}Note that if each function $f_{n}^{\ast}$ (with $n\geq1$) is
strictly convex on its domain, then the infimum in (\ref{r-1}), when finite,
is attained at at most one sequence $(x_{n}^{\ast})_{n\geq1}\subset
(\operatorname*{dom}f_{n}^{\ast})_{n\geq1}$ with $x^{\ast}=w^{\ast}$%
-$\sum_{n\geq1}x_{n}^{\ast}.$
\end{remark}

Because (\ref{r2}) is valid automatically for $x^{\ast}\in E^{\ast}%
\setminus\operatorname*{dom}f^{\ast}$, the problem is to see what is happening
for $x^{\ast}\in\operatorname*{dom}f^{\ast}\setminus\partial f\left(
\operatorname*{int}(\operatorname*{dom}f)\right)  .$

Taking $f_{k}=0$ for $k\geq n+1$ in Proposition \ref{p3} (v), its conclusion
is much weaker than the usual result mentioned at the beginning of this
section because nothing is said for sure for $x^{\ast}\in\operatorname*{dom}%
f^{\ast}\setminus\partial f\left(  \operatorname*{int}(\operatorname*{dom}%
f)\right)  .$

In the next two propositions we give complete descriptions for $f^{\ast}$ and
$g^{\ast},$ where $f$ and $g$ are provided in Proposition \ref{prop1} and
Example \ref{ex2}, respectively.

\begin{proposition}
\label{prop2}Let $f_{n}(x):=e^{\sigma_{n}x}$ for $x\in\mathbb{R}$ with
$0<\sigma_{n}\rightarrow\infty$, and $f=\sum_{n\geq1}f_{n}$. Assume that
$\operatorname*{dom}f\neq\emptyset$, and so $I:=(-\infty,-\alpha
)\subset\operatorname*{dom}f\subset\operatorname*{cl}I$ for some $\alpha
\in\mathbb{R}_{+}.$

\emph{(i)} Then $\partial f(\operatorname*{int}(\operatorname*{dom}%
f))=(0,\gamma)$, where $\gamma:=\sum_{n\geq1}\sigma_{n}e^{-\sigma_{n}\alpha
}\in\overline{\mathbb{R}}_{+},$ $\operatorname*{dom}f^{\ast}=\mathbb{R}_{+}$,
and
\begin{equation}
f^{\ast}(u)\leq\inf\bigg\{\sum_{n\geq1}f_{n}^{\ast}(u_{n})\mid(u_{n})_{n\geq
1}\subset\left(  \operatorname*{dom}f_{n}^{\ast}\right)  _{n\geq1}%
,\ u=\sum_{n\geq1}u_{n}\bigg\}<\infty\quad\forall u\in\mathbb{R}_{+}.
\label{r-70}%
\end{equation}

\emph{(ii)} Let $u\in\mathbb{R}_{+}$ $(=\operatorname*{dom}f^{\ast})$. Then
$u\in\lbrack0,\gamma]\cap\mathbb{R}$ if and only if
\begin{equation}
f^{\ast}(u)=\min\bigg\{\sum_{n\geq1}f_{n}^{\ast}(u_{n})\mid(u_{n})_{n\geq
1}\subset\left(  \operatorname*{dom}f_{n}^{\ast}\right)  _{n\geq1}%
,\ u=\sum_{n\geq1}u_{n}\bigg\}, \label{r7}%
\end{equation}
or, equivalently,%
\begin{equation}
f^{\ast}(u)=\min\bigg\{\sum_{n\geq1}u_{n}(\ln u_{n}-1)\mid(u_{n})_{n\geq
1}\subset\mathbb{R}_{+},\ u=\sum_{n\geq1}\sigma_{n}u_{n}\bigg\}. \label{r9}%
\end{equation}
More precisely, the minimum in (\ref{r9}) is attained for $u_{n}=0$ $(n\geq1)$
when $u=0$, for $u_{n}=e^{\sigma_{n}x}$ $(n\geq1)$ when $u=f^{\prime}(x)$ with
$x\in I$, for $u_{n}=e^{-\sigma_{n}\alpha}$ when $u=\gamma$ $(<\infty)$ (in
which case $\alpha\in\mathbb{R}_{+}^{\ast}$, $-\alpha\in\operatorname*{dom}%
\partial f=\operatorname*{dom}f$ and $f_{-}^{\prime}(-\alpha)=\gamma$).

\emph{(iii)} Let $u\in\mathbb{R}_{+}$ $(=\operatorname*{dom}f^{\ast})$. Then
\begin{align}
f^{\ast}(u)  &  =\inf\bigg\{\sum_{n\geq1}f_{n}^{\ast}(u_{n})\mid(u_{n}%
)_{n\geq1}\subset\left(  \operatorname*{dom}f_{n}^{\ast}\right)  _{n\geq
1},\ u=\sum_{n\geq1}u_{n}\bigg\}\label{r-13a}\\
&  =\inf\bigg\{\sum_{n\geq1}u_{n}(\ln u_{n}-1)\mid(u_{n})_{n\geq1}%
\subset\mathbb{R}_{+},\ u=\sum_{n\geq1}\sigma_{n}u_{n}\bigg\}. \label{r-13b}%
\end{align}

\end{proposition}

Proof. Let $\alpha\in\mathbb{R}_{+}$ be such that $I:=(-\infty,-\alpha
)=\operatorname*{int}(\operatorname*{dom}f)\subset\operatorname*{dom}%
f\subset\operatorname*{cl}I$ (see Proposition \ref{prop1}), and take
$\gamma:=\sum_{n\geq1}\sigma_{n}e^{-\sigma_{n}\alpha}\in(0,\infty]$.

(i) The equality $\partial f(\operatorname*{int}(\operatorname*{dom}%
f))=(0,\gamma)$ is proved in Proposition \ref{prop1} (ii). Because the
conjugate of the exponential function is given by
\begin{equation}
\exp^{\ast}(u)=\left\{
\begin{array}
[c]{ll}%
\infty & \text{if }u\in\mathbb{R}_{-}^{\ast},\\
u(\ln u-1) & \text{if }u\in\mathbb{R}_{+},
\end{array}
\right.  \label{r-exps}%
\end{equation}
where $0\cdot\ln0:=0$, we have $f_{n}^{\ast}(u)=\frac{u}{\sigma_{n}}(\ln
\frac{u}{\sigma_{n}}-1)$ for $u\in\mathbb{R}_{+}$ and $f_{n}^{\ast}(u)=\infty$
for $u\in\mathbb{R}_{-}^{\ast}.$

The first inequality in (\ref{r-70}) is given in (\ref{r-1}), while for the
second inequality just take $u_{1}:=u\in\mathbb{R}_{+}$ and $u_{n}:=0$ for
$n\geq2$. For $u\in\mathbb{R}_{-}^{\ast}$ we have that $f^{\ast}(u)=\sup
_{x\in\mathbb{R}}[xu-f(x)]\geq\lim_{x\rightarrow-\infty}[xu-f(x)]=\infty.$ It
follows that $\operatorname*{dom}f^{\ast}=\mathbb{R}_{+}.$

(ii) Consider $u\in\mathbb{R}_{+}$. Clearly, $f^{\ast}(0)=-\inf f=0$, and so
(\ref{r7}) and (\ref{r9}) hold in this case, with attainment for $u_{n}=0$
$(n\geq1)$. For $u\in(0,\gamma)=\partial f(\operatorname*{int}%
(\operatorname*{dom}f))=f^{\prime}(I)$, there exists $x\in I$ such that
$f^{\prime}(x)=u$, and so
\begin{align}
f^{\ast}(u)  &  =f^{\ast}(f^{\prime}(x))=\sup_{x\in\mathbb{R}}%
[xu-f(x)]=xf^{\prime}(x)-f(x)=\sum_{n\geq1}e^{\sigma_{n}x}\left(  \sigma
_{n}x-1\right) \label{r10a}\\
&  =\sum_{n\geq1}f_{n}^{\ast}(f_{n}^{\prime}(x))=\min\bigg\{\sum_{n\geq1}%
f_{n}^{\ast}(u_{n})\mid(u_{n})_{n\geq1}\subset\mathbb{R}_{+},\ u=\sum_{n\geq
1}u_{n}\bigg\}, \label{r10b}%
\end{align}
the last two equalities being given by Proposition \ref{p3} (v). Hence
(\ref{r7}) and (\ref{r9}) hold in this case, too, the attainment in (\ref{r9})
being for $u_{n}=e^{\sigma_{n}x}$ $(n\geq1)$. Hence, if $\gamma=\infty$ we
have that (\ref{r7}) (therefore, also (\ref{r9})) holds for all $u\in
\mathbb{R}_{+}=\operatorname*{dom}f^{\ast}=[0,\gamma]\cap\mathbb{R}.$

Assume that $\gamma<\infty$; then, by Proposition \ref{prop1}, we have that
$\alpha\in\mathbb{R}_{+}^{\ast}$ and $\gamma=f_{-}^{\prime}(-\alpha)$. Take
$u\geq\gamma$. Since $\psi^{\prime}(x)=u-f^{\prime}(x)>0$ for every $x\in I$,
where $\psi(x):=xu-f(x)$, it follows that $\psi$ is increasing on
$(-\infty,-\alpha]$, and so
\begin{equation}
f^{\ast}(u)=\sup_{x\in(-\infty,-\alpha]}\psi(x)=\psi(-\alpha)=-\alpha
u-f(-\alpha). \label{r-12}%
\end{equation}
Since $f^{\ast}$ is continuous on $\mathbb{R}_{+}^{\ast}=\operatorname*{int}%
(\operatorname*{dom}f^{\ast})$ and $\lim_{x\uparrow-\alpha}f^{\prime
}(x)=\gamma$, taking the limit for $x\uparrow-\alpha$ in (\ref{r10a}), we get
\[
f^{\ast}(\gamma)=-\alpha\gamma-f(-\alpha)=-\alpha\sum_{n\geq1}\sigma
_{n}e^{-\sigma_{n}\alpha}-\sum_{n\geq1}e^{-\sigma_{n}\alpha}=\sum_{n\geq
1}f_{n}^{\ast}(f_{n}^{\prime}(-\alpha)).
\]
Assume now that (\ref{r7}) holds for some $u>\gamma$, that is there exists
$(u_{n})_{n\geq1}\subset\mathbb{R}_{+}$ such that $u=\sum_{n\geq1}u_{n}$ and
$f^{\ast}(u)=\sum_{n\geq1}f_{n}^{\ast}(u_{n})$. Because $u\in\partial
f(-\alpha)=[\gamma,\infty)$, by Proposition \ref{prop1} (iv), we obtain that
$u_{n}\in\partial f_{n}(-\alpha)=\left\{  \sigma_{n}e^{-\sigma_{n}\alpha
}\right\}  ,$ whence $u=\sum_{n\geq1}\sigma_{n}e^{-\sigma_{n}\alpha}=\gamma$.
This contradiction proves that (\ref{r7}) holds if and only if $u\in
\lbrack0,\gamma]\cap\mathbb{R}.$

(iii) Because by (ii) the conclusion is clearly true for $u\in\lbrack
0,\gamma]\cap\mathbb{R}$, we may (and do) assume that $\gamma<\infty$. Take
$u>\gamma$ and denote by $F(u)$ the (real) number in the RHS of (\ref{r-13b}).
There exists $\overline{n}\geq1$ such that $\sigma_{n}\geq u$ for
$n\geq\overline{n}$, and fix $n\geq\overline{n}$. For $q\in\mathbb{N}^{\ast}$,
since $v:=u-\gamma+\sum_{k\geq n+1}\sigma_{k}e^{-\sigma_{k}\alpha}\in(0,u)$,
there exists a unique $\lambda_{q}\in\mathbb{R}_{+}^{\ast}$ such that
$\sum_{k=n+1}^{n+q}\sigma_{k}e^{-\sigma_{k}\lambda_{q}}=v$. It follows that
$\lambda_{q}<\lambda_{q+1}<\alpha$ for all $q\geq1;$ this follows easily by
contradiction. Therefore, $\lambda_{q}\uparrow\mu$ with $\mu\leq\alpha$.
Setting $u_{k}^{q}:=e^{-\sigma_{k}\alpha}$ for $k\in\overline{1,n},$
$u_{k}^{q}:=e^{-\sigma_{k}\lambda_{q}}$ for $k\in\overline{n+1,n+q}$ and
$u_{k}^{q}:=0$ for $k>n+q$, we have that $\sum_{k\geq1}\sigma_{k}u_{k}^{q}=u;$
taking into account (\ref{r-12}), we get%
\begin{align*}
f^{\ast}(u)  &  \leq F(u)\leq\sum_{k\geq1}u_{k}^{q}(\ln u_{k}^{q}%
-1)=\sum_{k=1}^{n}e^{-\sigma_{k}\alpha}(-\sigma_{k}\alpha-1)+\sum
_{k=n+1}^{n+q}e^{-\sigma_{k}\lambda_{q}}(-\sigma_{k}\lambda_{q}-1)\\
&  =(\lambda_{q}-\alpha)\sum_{k=1}^{n}\sigma_{k}e^{-\sigma_{k}\alpha}%
-\lambda_{q}u-\sum_{k=1}^{n}e^{-\sigma_{k}\alpha}-\sum_{k=n+1}^{n+q}%
e^{-\sigma_{k}\lambda_{q}}=f^{\ast}(u)+\Lambda_{q}^{n},
\end{align*}
where%
\begin{equation}
\Lambda_{q}^{n}:=(\alpha-\lambda_{q})\bigg(u-\sum_{k=1}^{n}\sigma
_{k}e^{-\sigma_{k}\alpha}\bigg)+\sum_{k\geq n+1}e^{-\sigma_{k}\alpha}%
-\sum_{k=n+1}^{n+q}e^{-\sigma_{k}\lambda_{q}}; \label{r-14}%
\end{equation}
hence $\Lambda_{q}^{n}\geq0$ for all $n\geq\overline{n}$ and $q\geq1.$

Assume that $\mu<\alpha$. Since $\lambda_{q}<\mu$, we have that $\sum
_{k=n+1}^{n+q}e^{-\sigma_{k}\lambda_{q}}\geq\sum_{k=n+1}^{n+q}e^{-\sigma
_{k}\mu}\rightarrow\infty$ for $q\rightarrow\infty$. From (\ref{r-14}) we get
the contradiction $0\leq\lim_{q\rightarrow\infty}\Lambda_{q}^{n}=-\infty$.
Hence $\mu=\alpha$. Using again (\ref{r-14}) we obtain that $\limsup
_{q\rightarrow\infty}\Lambda_{q}^{n}\leq\sum_{k\geq n+1}e^{-\sigma_{k}\alpha}%
$, and so $F(u)\leq f^{\ast}(u)+\sum_{k\geq n+1}e^{-\sigma_{k}\alpha}$ for
every $n\geq\overline{n}$. It follows that $F(u)\leq f^{\ast}(u).$ Therefore,
$F(u)=f^{\ast}(u)$, and so (\ref{r-13b}) holds. The proof is complete.
\hfill$\square$

\medskip

Observe that the condition $\sigma_{n}>0$ in Proposition \ref{prop2} is not
essential; because $\sigma_{n}\rightarrow\infty,$ $\sigma_{n}>0$ for some
$n_{0}\geq1$ and every $n\geq n_{0}$. Indeed, apply Proposition \ref{prop2}
for $g:=\sum_{n\geq n_{0}}f_{n}$, then the usual duality results for
$f=f_{1}+\ldots+f_{n_{0}-1}+g$. Of course, in the conclusion one must replace
$[0,\gamma]\cap\mathbb{R}$ by $\operatorname*{cl}[\partial
f(\operatorname*{int}(\operatorname*{dom}f))].$

One could ask, as for Theorem \ref{t-zt33}, if the condition
$\operatorname*{int}(\operatorname*{dom}f)\neq\emptyset$ is essential in
Proposition \ref{p3} (v). The next result proves that this is the case.

\begin{proposition}
\label{p4}Let $g,g_{n}$ be as in Example \ref{ex2}, where $(\varsigma
_{n})_{n\geq1}\subset\mathbb{R}_{+}$ is such that $\varsigma_{n}%
/n\rightarrow\infty$. Then $\operatorname*{dom}g^{\ast}=\mathbb{R}_{+}%
\times\mathbb{R}$ and
\begin{equation}
g^{\ast}(u,v)=\inf\bigg\{\sum_{n\geq1}\gamma_{n}(\ln\gamma_{n}-1)\mid
\gamma_{n}\geq0,\ u=\sum_{n\geq1}n\gamma_{n},\ v=\sum_{n\geq1}(-1)^{n}%
\varsigma_{n}\gamma_{n}\bigg\} \label{r8}%
\end{equation}
for all $(u,v)\in\left(  \mathbb{R}_{+}^{\ast}\times\mathbb{R}\right)
\cup\{(0,0)\}$, while for $u=0\neq v$ the term in the RHS of (\ref{r8}) is
$\infty$. Moreover, for $u=0=v$ the infimum in (\ref{r8}) is attained, while
for $u\in\mathbb{R}_{+}^{\ast}$ the infimum in (\ref{r8}) is attained if and
only if $\overline{v}:=\sum_{n\geq1}(-1)^{n}\varsigma_{n}\left(
\frac{1+2u-\sqrt{4u+1}}{2u}\right)  ^{n}\in\mathbb{R}$ and $v=\overline{v}.$
\end{proposition}

Proof. On one hand, because $g(x,y)=f(x)+\iota_{\{0\}}(y)$, where $f$ is
defined in (\ref{r-2}), we have that $g^{\ast}(u,v)=f^{\ast}(u)$ for
$(u,v)\in\mathbb{R}^{2}$. On the other hand $g_{n}^{\ast}(u,v)=\frac{u}{n}%
(\ln\frac{u}{n}-1)$ for $u\in\mathbb{R}_{+}$ and $v=(-1)^{n}\varsigma_{n}/n,$
while $g_{n}^{\ast}(u,v)=\infty$ otherwise. By Proposition \ref{p3} (i) we
have that
\begin{equation}
g^{\ast}(u,v)\leq G(u,v):=\inf\bigg\{\sum_{n\geq1}\gamma_{n}(\ln\gamma
_{n}-1)\mid\gamma_{n}\geq0,\ u=\sum_{n\geq1}n\gamma_{n},\ v=\sum_{n\geq
1}(-1)^{n}\varsigma_{n}\gamma_{n}\bigg\} \label{r6}%
\end{equation}
for all $(u,v)\in\mathbb{R}_{+}\times\mathbb{R}.$

It is clear that for $u=v=0$ one has equality with attained infimum (for
$\gamma_{n}=0$ for every $n\geq1$), while for $u=0\neq v$ the RHS term of
(\ref{r6}) is $\infty.$

Applying Proposition \ref{prop2} for $\sigma_{n}:=n$ $(n\geq1)$ we have that
$\alpha=0$ and $\gamma=\infty$; moreover, for $u\in\mathbb{R}_{+}^{\ast},$
\[
f^{\ast}(u)=\min\bigg\{\sum_{n\geq1}\gamma_{n}(\ln\gamma_{n}-1)\mid\gamma
_{n}\geq0,\ u=\sum_{n\geq1}n\gamma_{n}\bigg\},
\]
which is attained only for the sequence $(\overline{\gamma}_{n}^{u})_{n\geq1}%
$, where
\begin{equation}
\overline{\gamma}_{n}^{u}:=\bigg( \frac{1+2u-\sqrt{4u+1}}{2u}\bigg) ^{n}%
\quad(n\geq1). \label{r-gbnu}%
\end{equation}
Consequently, we have equality in (\ref{r6}) with attained infimum if and only
if $\overline{v}\in\mathbb{R}$ and $v=\overline{v}.$

Let $(u,v)\in\mathbb{R}_{+}^{\ast}\times\mathbb{R}$ and fix $\varepsilon>0.$
Then there exists some $\overline{n}\geq1$ such that $\sum_{k=1}^{\overline
{n}}\bar{\gamma}_{k}^{u}(\ln\bar{\gamma}_{k}^{u}-1)<f^{\ast}(u)+\varepsilon
/2$. Set $\overline{u}:=\sum_{k=1}^{\overline{n}}k\bar{\gamma}_{k}^{u}$ and
$\overline{v}:=\sum_{k=1}^{\overline{n}}(-1)^{k}\varsigma_{k}\bar{\gamma}%
_{k}^{u}$; then $u^{\prime}:=u-\overline{u}>0$ and $v^{\prime}:=v-\overline
{v}\in\mathbb{R}$. Observe that for a fixed $n\in\mathbb{N}^{\ast}$, since
$g_{n}$ is finite (and continuous) one has $(g_{n}+g_{n+1})^{\ast}=g_{n}%
^{\ast}\square g_{n+1}^{\ast}$ with exact convolution, and so
\[
\operatorname*{dom}(g_{n}+g_{n+1})^{\ast}=\operatorname*{dom}g_{n}^{\ast
}+\operatorname*{dom}g_{n+1}^{\ast}=\mathbb{R}_{+}(n,(-1)^{n}\varsigma
_{n})+\mathbb{R}_{+}(n+1,(-1)^{n+1}\varsigma_{n+1}).
\]
Take $\overline{n}\geq1$ such that $u^{\prime}\varsigma_{n}\geq n\left\vert
v^{\prime}\right\vert $ for $n\geq\overline{n}$. Then
\[
\gamma_{n}^{\prime}:=\frac{u^{\prime}\varsigma_{n+1}-(-1)^{n+1}(n+1)v^{\prime
}}{n\varsigma_{n+1}+(n+1)\varsigma_{n}}\geq0,\quad\gamma_{n+1}^{\prime}%
:=\frac{u^{\prime}\varsigma_{n}-(-1)^{n}nv^{\prime}}{n\varsigma_{n+1}%
+(n+1)\varsigma_{n}}\geq0,
\]
and $(u^{\prime},v^{\prime})=\gamma_{n}^{\prime}(n,(-1)^{n}\varsigma
_{n})+\gamma_{n+1}^{\prime}(n+1,(-1)^{n+1}\varsigma_{n+1})\in
\operatorname*{dom}(g_{n}+g_{n+1})^{\ast}$. Moreover,
\begin{align*}
(g_{n}+g_{n+1})^{\ast}(u^{\prime},v^{\prime})  &  =g_{n}^{\ast}(n\gamma
_{n}^{\prime},(-1)^{n}\varsigma_{n}\gamma_{n}^{\prime})+g_{n+1}^{\ast
}((n+1)\gamma_{n+1}^{\prime},(-1)^{n+1}\varsigma_{n+1}\gamma_{n+1}^{\prime})\\
&  =\gamma_{n}^{\prime}(\ln\gamma_{n}^{\prime}-1)+\gamma_{n+1}^{\prime}%
(\ln\gamma_{n+1}^{\prime}-1).
\end{align*}
Since for $n\geq\overline{n}$ one has $0\leq\gamma_{n}^{\prime}\leq
\frac{u^{\prime}\varsigma_{n+1}+(n+1)\left\vert v^{\prime}\right\vert
}{n\varsigma_{n+1}}\leq u^{\prime}\frac{1}{n}+\left\vert v^{\prime}\right\vert
\frac{2}{\varsigma_{n+1}}\rightarrow0$ and $0\leq\gamma_{n+1}^{\prime}%
\leq\frac{u^{\prime}\varsigma_{n}+n\left\vert v^{\prime}\right\vert
}{(n+1)\varsigma_{n}}=u^{\prime}\frac{1}{n+1}+\left\vert v^{\prime}\right\vert
\frac{1}{\varsigma_{n}}\rightarrow0$, there exists $m\geq\overline{n}$ such
that $\gamma_{m}^{\prime}(\ln\gamma_{m}^{\prime}-1)+\gamma_{m+1}^{\prime}%
(\ln\gamma_{m+1}^{\prime}-1)<\varepsilon/2$; set $\gamma_{k}:=\bar{\gamma}%
_{k}^{u}$ for $k\in\overline{1,\overline{n}}$, $\gamma_{k}:=\gamma_{k}%
^{\prime}$ for $k\in\{m,m+1\}$ and $\gamma_{k}:=0$ otherwise. Then $\gamma
_{n}\geq0$ for all $n\geq1$, $u=\sum_{n\geq1}n\gamma_{n}$ and $v=\sum_{n\geq
1}(-1)^{n}\varsigma_{n}\gamma_{n}$; moreover,%
\begin{align*}
G(u,v)  &  \leq\sum_{n\geq1}\gamma_{n}(\ln\gamma_{n}-1)=\sum_{k=1}%
^{\overline{n}}\bar{\gamma}_{k}^{u}(\ln\bar{\gamma}_{k}^{u}-1)+\gamma
_{m}^{\prime}(\ln\gamma_{m}^{\prime}-1)+\gamma_{m+1}^{\prime}(\ln\gamma
_{m+1}^{\prime}-1)\\
&  \leq f^{\ast}(u)+\varepsilon/2+\varepsilon/2=f^{\ast}(u)+\varepsilon
=g^{\ast}(u,v)+\varepsilon.
\end{align*}
It follows that $G(u,v)\leq g^{\ast}(u,v)$, and so (\ref{r8}) holds. The proof
is complete. \hfill$\square$

\begin{remark}
Observe that, depending on $(\varsigma_{n})_{n\geq1}$, the set of those $u>0$
for which the infimum in the RHS of (\ref{r8}) is attained for some
$v\in\mathbb{R}$ could be $\mathbb{R}_{+}^{\ast},$ a proper subset of
$\mathbb{R}_{+}^{\ast}$, or the empty set. For example, when $\varsigma
_{n}=n^{k}$ with $k>1$ the series $\sum_{n\geq1}(-1)^{n}\varsigma_{n}%
\overline{\gamma}_{n}^{u}$ is convergent. If $\varsigma_{n}=e^{n\alpha}$ with
$\alpha>0$ the series $\sum_{n\geq1}(-1)^{n}\varsigma_{n}\overline{\gamma}%
_{n}^{u}$ is convergent iff $\alpha<\ln\frac{1+2u+\sqrt{4u+1}}{2u}$. If
$\varsigma_{n}=e^{n^{2}}$, the series $\sum_{n\geq1}(-1)^{n}\varsigma
_{n}\overline{\gamma}_{n}^{u}$ is not convergent. (Here $\overline{\gamma}%
_{n}^{u}$ is defined in (\ref{r-gbnu}).)
\end{remark}

Having in view Propositions \ref{prop2} and \ref{p4}, the following question
is natural:

\begin{question}
\label{q2}Let $f,f_{n}\in\Lambda(\mathbb{R}^{p})$ with $p\in\mathbb{N}^{\ast}$
and $f(x)=\sum_{n\geq1}f_{n}(x)$ for every $x\in\mathbb{R}^{p}$; is it true
that
\[
f^{\ast}(u)=\inf\bigg\{\sum_{n\geq1}f_{n}^{\ast}(u_{n})\mid(u_{n})_{n\geq
1}\subset(\operatorname*{dom}f_{n}^{\ast})_{n\geq1},\ u=\sum_{n\geq1}%
u_{n}\bigg\}
\]
for all $u\in\operatorname*{int}(\operatorname*{dom}f^{\ast})?$
\end{question}

M. Valadier in \cite{Val:70} defines the continuous inf-convolution for a
family $(f_{t})_{t\in T}$ of proper lower semi\-continuous functions defined
on $\mathbb{R}^{p}$, where $(T,\mathcal{T},\mathcal{\mu)}$ is a measure space
with $\mu\geq0$ being $\sigma$-finite. In the case in which $T:=\mathbb{N}%
^{\ast}$, $\mathcal{T}:=2^{\mathbb{N}}$ and $\mu:\mathcal{T}\rightarrow
\lbrack0,\infty]$ is defined by $\mu(A):=\operatorname*{card}A,$
\cite[Th.\ 7]{Val:70} has the following (equivalent) statement:

\medskip

\textbf{Theorem} \emph{Let $(g_{n})_{n\geq1}\subset\Gamma(\mathbb{R}^{p})$ be
such that the series $\sum_{n\geq1}\left\vert g_{n}^{\ast}(u)\right\vert
<\infty$ for every $u\in\mathbb{R}^{p}.$ Then the function $\widetilde
{g}:\mathbb{R}^{p}\rightarrow\overline{\mathbb{R}},$ defined by
\[
\widetilde{g}(x):=\inf\bigg\{
\sum_{n\geq1}g_{n}(x_{n})\mid(x_{n})_{n\geq1}\subset\mathbb{R}^{p}%
,\ \sum_{n\geq1}\left\Vert x_{n}\right\Vert <\infty,\ x=\sum_{n\geq1}%
x_{n}\bigg\}  ,
\]
belongs to $\Gamma(\mathbb{R}^{p}),$ the infimum above is attained for every
$x\in\mathbb{R}^{p}$, and $\widetilde{g}^{\ast}=\sum_{n\geq1}g_{n}^{\ast}$.}

\medskip

Taking $(f_{n})_{n\geq1}\subset\Gamma(\mathbb{R}^{p})$, and setting
$g_{n}:=f_{n} ^{\ast}$ (hence $g_{n}^{\ast}=f_{n}$), the hypothesis of
\cite[Th.\ 7]{Val:70} implies that $\operatorname*{dom}f=\mathbb{R}^{p}$ (and
of course $f,f_{n}$ are continuous on $\mathbb{R}^{p}$), an hypothesis which
is stronger than that of Proposition \ref{p3} (v), but also the conclusion of
\cite[Th.\ 7]{Val:70} is stronger. Clearly, \cite[Th.\ 7]{Val:70} (above) can
not be applied in the previous examples (because $f$ is not finite-valued), as
well as for the example in the next section.\footnote{We thank Prof.\ L.
Thibault for bringing to our attention the reference \cite{Val:70}.}

\section{An application to entropy minimization}

As mentioned in Introduction, in Statistical Physics (Statistical Mechanics)
one has to minimize $\sum_{i\in I}n_{i}(\ln n_{i}-1)$ with the constraints
$\sum_{i\in I}n_{i}=N$ and $\sum_{i\in I}n_{i}\varepsilon_{i}=\varepsilon$,
where $I$ is a countable set and $n_{i}$ are nonnegative integers. In this
context consider $h_{i}:=\exp\circ A_{i}$, where $A_{i}:\mathbb{R}%
^{2}\rightarrow\mathbb{R}$ is defined by $A_{i}(x,y):=x+\varepsilon_{i}y$, and
$h:=\sum_{i\in I}h_{i}.$ Because $h_{i}>0$, we have that $\sum_{i\in I}%
h_{i}=\sum_{j\in J}h_{p(j)}$ for every bijection $p:J\rightarrow I$. So, we
(can) take $I=\mathbb{N}^{\ast}$, the set of positive integers, $(\sigma
_{n})_{n\geq1}\subset\mathbb{R},$ $A_{n}(x,y):=x+\sigma_{n}y$, $h_{n}%
:=\exp\circ A_{n}$ and $h:=\sum_{n\geq1}h_{n}$. Hence
\[
h(x,y)=\sum_{n\geq1}e^{x+\sigma_{n}y}=e^{x}\sum_{n\geq1}e^{\sigma_{n}y}%
=e^{x}f(y)
\]
with $f:\mathbb{R}\rightarrow\overline{\mathbb{R}}$, $f(y):=\sum_{n\geq
1}e^{\sigma_{n}y}$. By Proposition \ref{prop1}, when $\operatorname*{dom}%
h=\mathbb{R}\times\operatorname*{dom}f\neq\emptyset$, we may (and do) assume
that $\sigma_{n}\rightarrow\infty$, in which case there exists $\alpha
\in\mathbb{R}_{+}$ such that $I:=(-\infty,-\alpha)=\operatorname*{int}%
(\operatorname*{dom}f)\subset\operatorname*{dom}f\subset\operatorname*{cl}I.$
Since $h_{n}$ is differentiable on $\mathbb{R}^{2}$ [with $\nabla
h_{n}(x,y)=e^{x+\sigma_{n}y}(1,\sigma_{n})$], by Theorem \ref{t-zt33} (and
Corollary \ref{c-zc31}), $h$ is differentiable on $\operatorname*{int}%
(\operatorname*{dom}h)=\mathbb{R}\times I$ and
\[
\nabla h(x,y)=\sum_{n\geq1}e^{x+\sigma_{n}y}(1,\sigma_{n})=\bigg(\sum_{n\geq
1}e^{x+\sigma_{n}y},\sum_{n\geq1}\sigma_{n}e^{x+\sigma_{n}y}\bigg)=\left(
e^{x}f(y),e^{x}f^{\prime}(y)\right)
\]
for all $(x,y)\in\mathbb{R}\times I$. Moreover, because $\operatorname{Im}%
A_{n}=\mathbb{R}$ and $A_{n}^{\ast}w=w(1,\sigma_{n})$ for $w\in\mathbb{R}$,
\[
h_{n}^{\ast}(u,v)=\min\left\{  \exp^{\ast}(w)\mid A_{n}^{\ast}w=(u,v)\right\}
=\left\{
\begin{array}
[c]{ll}%
u(\ln u-1) & \text{if }u\geq0\text{ and }v=u\sigma_{n},\\
\infty & \text{otherwise.}%
\end{array}
\right.
\]

Using Proposition \ref{p3} (v), it follows that for $(u,v)=\nabla f(x,y)$ with
$(x,y)\in\operatorname*{int}(\operatorname*{dom}f),$
\begin{equation}
h^{\ast}(u,v)=\min\bigg\{\sum_{n\geq1}u_{n}(\ln u_{n}-1)\mid(u_{n})_{n\geq
1}\in\mathbb{R}_{+},\ \sum_{n\geq1}u_{n}=u,\ \sum_{n\geq1}u_{n}\sigma
_{n}=v\bigg\} \label{r4}%
\end{equation}
for every $(u,v)\in\partial h\left(  \operatorname*{int}(\operatorname*{dom}%
f)\right)  $; moreover, because $h_{n}^{\ast}$ is strictly convex, for
$(u,v)\in\partial h(x,y)$ with $(x,y)\in\operatorname*{int}%
(\operatorname*{dom}h)$, the minimum in (\ref{r4}) is realized at the unique
sequence $(\overline{u}_{n})_{n\geq1}=(e^{x+\sigma_{n}y})_{n\geq1}.$

We apply the preceding considerations for the following example taken from
\cite[p.\ 10]{PatBea:11}:

\medskip

\textquotedblleft\ $\varepsilon(n_{x},n_{y},n_{z})=\frac{h^{2}}{8mL^{2}}%
(n_{x}^{2}+n_{y}^{2}+n_{z}^{2});\quad n_{x},n_{y},n_{z}=1,2,3,\ldots$ (5)

\medskip

\noindent where $h$ is Planck's constant and $m$ the mass of the
particle,\textquotedblright\ and $L$ is the side of a cubical box.

Hence, $\sigma_{k,l,m}:=\kappa(k^{2}+l^{2}+m^{2})$ with $k,l,m\in
\mathbb{N}^{\ast}$. We take $\kappa=1$ to (slightly) simplify the calculation.
It follows that
\[
h(x,y)=\sum_{k,l,m\geq1}e^{x+y(k^{2}+l^{2}+m^{2})}=e^{x}\sum_{k,l,m\geq
1}e^{yk^{2}}e^{yl^{2}}e^{ym^{2}}=e^{x}\bigg(\sum_{k\geq1}e^{yk^{2}}%
\bigg)^{3}.
\]

Clearly, $\operatorname*{dom}h=\mathbb{R}\times\mathbb{R}_{-}^{\ast
}=\operatorname*{int}(\operatorname*{dom}h)$. Let us consider
\[
f:\mathbb{R}\rightarrow\overline{\mathbb{R}},\quad f(y):={\sum}_{n\geq
1}e^{n^{2}y};
\]
hence $h(x,y)=e^{x}\left[  f(y)\right]  ^{3}$. As observed in Example
\ref{ex1}, $I:=\operatorname*{dom}f=\mathbb{R}_{-}^{\ast}$, and so
$\operatorname*{dom}h=\mathbb{R}\times I$. The series $\sum_{k\geq1}e^{yk^{2}%
}$, as well as the series $\sum_{k\geq1}(k^{2})^{p}e^{yk^{2}}$ with
$p\in\mathbb{N}^{\ast}$, are uniformly convergent on the interval
$(-\infty,-\gamma]$ for every $\gamma>0$ (because $0\leq e^{yk^{2}}\leq
e^{-\gamma k^{2}}$ for every $y\in(-\infty,-\gamma]$). It follows that
$f^{(p)}(y)=\sum_{k\geq0}(k^{2})^{p}e^{yk^{2}}$ for every $p\in\mathbb{N}$
$($with $f^{(0)}:=f,$ $f^{(1)}:=f^{\prime}$ ...) and $\lim_{y\rightarrow
-\infty}f^{(p)}(y)=0$ for $p\in\mathbb{N}$. Moreover, $\lim_{y\rightarrow
0-}f(y)=\infty$. This is because $f(y)\geq\sum_{k=1}^{n}e^{yk^{2}}$ for every
$n\geq1$ and $\lim_{y\rightarrow0-}\sum_{k=1}^{n}e^{yk^{2}}=n$. Hence
$h|_{\operatorname*{dom}h}\in C^{\infty}(\operatorname*{dom}h)$. We know that
$h$ is strictly convex on its domain (as the sum of a series of strictly
convex functions). In fact $\ln f$ (with $\ln\infty:=\infty$) is proper,
convex and lsc; $\ln f$ is even strictly convex on $\operatorname*{dom}f=I$.
Indeed, $(\ln f)^{\prime}=f^{\prime}/f>0$ and $(\ln f)^{\prime\prime
}=(f^{\prime\prime}f-(f^{\prime})^{2})/f^{2}>0$ on $I$; just use Schwarz
inequality in $\ell^{2}$. Moreover, $\lim_{y\rightarrow-\infty}f^{\prime
}(y)/f(y)=1$ and $\eta:=\lim_{y\rightarrow0-}f^{\prime}(y)/f(y)=\infty$. The
first limit is (almost) obvious. The second limit exists because $(\ln
f)^{\prime}$ is increasing on $I$. In fact, for fixed $n>1$, we have that
$f^{\prime}(y)=\sum_{k\geq1}k^{2}e^{k^{2}y}\geq n^{2}f(y)-n^{2}\sum_{k=1}%
^{n}k^{2}e^{k^{2}y}$, and so $f^{\prime}(y)/f(y)\geq n^{2}-n^{2}\left(
\sum_{k=1}^{n}k^{2}e^{k^{2}y}\right)  /f(y)$ for $y<0$. Since $\lim
_{y\rightarrow0-}f(y)=\infty$, it follows that $\eta\geq n^{2}-n^{2}\left(
\sum_{k=1}^{n}k^{2}\right)  /\infty=n^{2}$. Therefore, $\eta=\infty$. It
follows that $\varphi:=f^{\prime}/f:\mathbb{R}_{-}^{\ast}\rightarrow
(1,\infty)$ is a bijection.

Let us first determine the conjugate of $\ln f$ which will be needed to
express the conjugate $h^{\ast}$ of $h$. Since the equation $(\frac{d}%
{dy}[vy-\ln f(y)]=)$ $v-\varphi(y)=0$ has the (unique) solution $y=\varphi
^{-1}(v)\in I$ for $v>1$, we obtain that
\[
(\ln f)^{\ast}(v)=\sup\{vy-\ln f(y)\mid y\in I\}=v\varphi^{-1}(v)-\ln\left[
f(\varphi^{-1}(v))\right]  \quad\forall v>1.
\]
Because $(\ln f)^{\ast}$ is lsc, we have that
\begin{align*}
(\ln f)^{\ast}(1)  &  =\lim_{v\rightarrow1+}(\ln f)^{\ast}(v)=\lim
_{y\rightarrow-\infty}\left[  y\varphi(y)-\ln f(y)\right] \\
&  =\lim_{y\rightarrow-\infty}\bigg[y\frac{1+\sum_{n\geq2}n^{2}e^{(n^{2}-1)y}%
}{1+\sum_{n\geq2}e^{(n^{2}-1)y}}-y-\ln\bigg(1+\sum_{n\geq2}e^{(n^{2}%
-1)y}\bigg)\bigg]=0.
\end{align*}
Finally,
\begin{align*}
(\ln f)^{\ast}(v)  &  =\sup_{y\in I}[vy-\ln f(y)]=\lim_{y\rightarrow-\infty
}[vy-\ln f(y)]=(\lim_{y\rightarrow-\infty}y)\left(  v-\lim_{y\rightarrow
-\infty}\frac{\ln f(y)}{y}\right) \\
&  =(-\infty)\left(  v-\lim_{y\rightarrow-\infty}\frac{f^{\prime}(y)}%
{f(y)}\right)  =(-\infty)(v-1)=\infty\quad\forall v<1.
\end{align*}

Let us determine now the conjugate of $h$ for $(u,v)\in\mathbb{R}^{2}$. Since
$f(y)\in\mathbb{R}_{+}^{\ast}$ for $y\in\mathbb{R}_{-}^{\ast},$
\begin{align*}
h^{\ast}(u,v)  &  =\sup_{(x,y)\in\operatorname*{dom}h}[xu+yv-h(x,y)]=\sup
_{y\in\mathbb{R}_{-}^{\ast}}\left(  yv+\sup_{x\in\mathbb{R}}\left[
xu-e^{x}[f(y)]^{3}\right]  \right) \\
&  =\sup_{y\in\mathbb{R}_{-}^{\ast}}\left(  yv+[f(y)]^{3}\exp^{\ast}\left(
\frac{u}{[f(y)]^{3}}\right)  \right)  .
\end{align*}
Hence $h^{\ast}(u,v)=\infty$ for $u\in\mathbb{R}_{-}^{\ast}$ and $h^{\ast
}(0,v)=\iota_{\mathbb{R}_{-}^{\ast}}^{\ast}(v)=\iota_{\mathbb{R}_{+}}(v)$. For
$u\in\mathbb{R}_{+}^{\ast}$ we have that
\begin{align*}
h^{\ast}(u,v)  &  =\sup_{y\in\mathbb{R}_{-}^{\ast}}\left(  yv+u(\ln
u-3\ln[f(y)]-1\right)  =u(\ln u-1)+\sup_{x\in\mathbb{R}}\left(  yv-3u\ln
[f(y)]\right) \\
&  =u(\ln u-1)+3u(\ln f)^{\ast}\left(  \frac{v}{3u}\right)  .
\end{align*}
In conclusion,
\[
h^{\ast}(u,v)=\left\{
\begin{array}
[c]{ll}%
u(\ln u-1)+3u(\ln f)^{\ast}\left(  \frac{v}{3u}\right)  & \text{if }%
v\geq3u>0,\\
0 & \text{if }u=0\leq v,\\
\infty & \text{if }u<0,\text{ or }v<0,\text{ or }0\leq v<3u.
\end{array}
\right.
\]
It follows that
\[
\left\{  (u,v)\in\mathbb{R}_{+}^{\ast}\times\mathbb{R}_{+}^{\ast}\mid
v>3u\right\}  =\operatorname*{int}(\operatorname*{dom}h^{\ast})\subset
\operatorname*{dom}h^{\ast}=\left\{  (u,v)\in\mathbb{R}_{+}\times
\mathbb{R}_{+}\mid v\geq3u\right\}  .
\]
Moreover, because $\nabla h(x,y)=e^{x}\left[  f(y)\right]  ^{3}\cdot
(1,3\varphi(y))$, we get
\[
\partial h(\operatorname*{int}(\operatorname*{dom}h))=\nabla h(\mathbb{R}%
\times I)=\operatorname*{int}(\operatorname*{dom}h^{\ast}).
\]

For $(u,v)\in\mathbb{R}^{2}$ consider the set%
\[
S(u,v):=\bigg\{(u_{k,l,m})_{k,l,m\geq1}\subset\mathbb{R}_{+}\mid
\sum_{k,l,m\geq1}u_{k,l,m}=u,~\sum_{k,l,m\geq1}(l^{2}+k^{2}+m^{2}%
)u_{k,l,m}=v\bigg\};
\]
clearly, $S(u,v)=\emptyset$ for $(u,v)\notin\operatorname*{dom}h^{\ast},$
$S(u,3u)=\{(u_{k,l,m})_{k,l,m\geq1}\subset\mathbb{R}_{+}\mid u_{1,1,1}:=u,$
$u_{k,l,m}:=0$ otherwise$\}$ for $u\geq0$ and $S(0,v)=\emptyset$ for $v>0.$

Applying Proposition \ref{p3} (v), for $(u,v)\in\operatorname*{int}%
(\operatorname*{dom}h^{\ast})=\partial h(\operatorname*{int}%
(\operatorname*{dom}h))$ we have that
\begin{equation}
h^{\ast}(u,v)=\min\bigg\{\sum_{k,l,m\geq1}u_{k,l,m}(\ln u_{k,l,m}%
-1)\mid(u_{k,l,m})_{k,l,m\geq1}\in S(u,v)\bigg\}, \label{r-exmb}%
\end{equation}
the minimum being attained uniquely for $(\overline{u}_{k,l,m}):=\frac
{u}{\left[  f(y)\right]  ^{3}}e^{(k^{2}+l^{2}+m^{2})y}$ $(k,l,m\geq1)$, where
$y<0$ is the solution of the equation $v/u=3f^{\prime}(y)/f(y)$; hence
$S(u,v)\neq\emptyset$ in this case. In the case $v=3u\geq0$, as seen above,
$S(u,3u)$ is a singleton and relation (\ref{r-exmb}) holds, too; for $v>0$
$(=u)$ we have that $S(0,v)=\emptyset\ $and $h^{\ast}(0,v)=0$.

Observe that the solution for the case $v=3u\geq0$ is not obtained by using
the (formal) method of Lagrange multipliers. Also note that even the solution
for $(u,v)\in\partial h(\operatorname*{int}(\operatorname*{dom}h))$ can not be
obtained from the results of J. M. Borwein and his collaborators because, even
if $\ell^{p}$-spaces can be regarded as $L^{p}(\Omega)$-spaces, the measure of
$\Omega$ is not finite (and, even more, the corresponding linear operators are
not continuous).

\medskip

\textbf{Acknowledgement.} C.Z. wishes to thank his co-author for introducing
him in entropy minimization problems related to Statistical Physics. C.Z. also
benefited from discussions with C.V. on the use of Convex Analysis in
Mechanics. Related to this, C.V. recalled the influence J. J. Moreau had on
his vision on Mechanics. J. J. Moreau was a member of the jury of C. V.'s PhD
thesis \cite{Val:73} defended at Universit\'{e} de Poitiers in 1973; they were
colleagues for a while, and remained good friends. Unfortunately, C.V. died in
November 2014.

Thanks go to Prof.\ M. Durea and Prof.\ D. Fortun\'{e} for their remarks on a
previous version of the manuscript.

\end{document}